\magnification=1200

\def\restriction{|}
\def\Vdash{{|\!\vdash}}
\let\frak\bf

\def\on{{\restriction}}
\def\forces{\Vdash}
\def\forcess{\mathop{\Vdash}}


\def\meager{{\cal M}}  
\def\measure0{{\cal N}}
\def\SMZ/{\hbox{$\Unif(\smz)$}}
\def\smz{{\cal S}}

\def\Add{{\bf Add}}
\def\Cov{{\bf Cov}}
\def\Unif{{\bf Unif}}

\parindent=0cm
\hfuzz=6pt

\def\begindent{\par\smallskip\begingroup\advance\parindent by
1cm\advance\rightskip by 1.5 cm plus 1 cm minus 1 cm\relax}
\def\ite#1 {\item{(#1)} }
\def\endent{\smallskip\endgroup}

\def\R{\hbox{I}\!\hbox{R}} 
\def\P{{\frak P}}

\def\Q{{\bf Q}}   
\def\U{{\cal U}} 

\def\b{{\frak b}}
\def\d{{\frak d}}
\def\J{{\cal J}}  
\def\c{{\frak c}}
\def\H{{\cal H}}  
\def\A{{\cal A}}  

\newcount\secno
\newcount\theono   
\newcount\chapno 
\def\smzDef/{0.1} 
\def\smzDefRem/{0.2} 
\def\defdef/{0.3} 
\def\defdefFact/{0.4} 
\def\GluzinDef/{0.5} 
\def\GluzinFact/{0.6} 
\def\onetheo/{0.7} 
\def\twotheo/{0.8} 
\def\SMZdef/{0.9} 
\def\RothsEq/{0.10} 
\def\rothTheo/{0.11} 
\def\rothFact/{0.12} 
\def\rothTproof/{0.13} 
\def\smallsmz/{0.14} 
\def\smallprime/{0.15} 
\def\bigsmz/{0.16} 
\def\sNotation/{0.17} 
\def\mNotation/{0.18} 
\def\emNot/{0.19} 
\def\ppDef/{1.1} 
\def\ppGame/{1.2} 
\def\pppDef/{1.3} 
\def\RPSdef/{1.4} 
\def\MillersLemma/{1.5} 
\def\ccFact/{1.6} 
\def\nonew/{1.7} 
\def\ccLemma/{1.8} 
\def\preservePP/{1.9} 
\def\oobDef/{1.10} 
\def\preserveoo/{1.11} 
\def\noCfact/{1.12} 
\def\strongoo/{1.13} 
\def\strongooIt/{1.14} 
\def\strooFact/{1.15} 
\def\Random/{1.16} 
\def\splitstat/{1.17} 
\def\rzeroone/{1.18} 
\def\HPerfDef/{2.1} 
\def\splitDef/{2.2} 
\def\fPerfRem/{2.3} 
\def\RPremark/{2.4} 
\def\splitFac/{2.5} 
\def\leDef/{2.6} 
\def\sameStem/{2.7} 
\def\lenFact/{2.8} 
\def\fusionFact/{2.9} 
\def\trimFact/{2.10} 
\def\joinFact/{2.11} 
\def\finSplit/{2.12} 
\def\SacksFact/{2.13} 
\def\PTcor/{2.14} 
\def\RPprop/{2.15} 
\def\RPproper/{2.16} 
\def\ptgenDef/{2.17} 
\def\ptgenFact/{2.18} 
\def\againRoths/{2.19} 
\def\PTinterpret/{2.20} 
\def\PTintFact/{2.21} 
\def\ultraLemma/{2.22} 
\def\ptwFact/{2.23} 
\def\ptwtwo/{2.24} 
\def\splitsDef/{2.25} 
\def\splitRem/{2.26} 
\def\splitOne/{2.27} 
\def\splitAll/{2.28} 
\def\funny/{2.29} 
\def\noincoo/{2.30} 
\def\ooproblem/{2.31} 
\def\hstarfact/{2.32} 
\def\reflect/{3.1} 
\def\smzisclub/{3.2} 
\def\addcor/{3.3} 
\def\smallproof/{3.4} 
\def\noincLemma/{3.5} 
\def\landDef/{3.6} 
\def\landRdef/{3.7} 
\def\noincP/{3.8} 
\def\bigproof/{3.9} 
\def \IteratedForcing {1}
\def \BlassShelah {2}
\def \Corazza {3}
\def \ShelahsPreservation {4}
\def \JSW {5}
\def \smzRapid {6}
\def \Kunen {7}
\def \Mapping {8}
\def \Rational {9}
\def \SomeProperties {10}
\def \Pawlikowski {11}
\def \ProprieteC {12}
\def \EigenschaftC {13}
\def \ProperForcing {14}
\def \ProperImproper {15}
\def \IteratedCohen {16}
\def \RealValued {17}

\catcode`@=11

\def\write#1#2{}
\newif\ifproofmode
\proofmodefalse            
\def\@nofirst#1{}
\def\neusection{\advance\secno by 1\relax \theono=0\relax}
\def\neuchap{\advance\chapno by 1\relax\secno=0\relax\theono=0\relax}
\neuchap
\def\labelit#1{\global\advance\theono by 1%
\global\edef#1/{%
\number\secno.\number\theono}%
}
\def\ppro#1/#2:{\labelit{#1}%
\smallbreak\noindent%
\@markit{#1}%
{\bf #2:}}
\def\@definition#1{\string\def\string#1{#1}
\expandafter\@nofirst\string\%
(\the\pageno)}

\def\@markit#1{\leavevmode
\ifproofmode\llap{{\tt \expandafter\@nofirst\string#1\ }}\fi
{\bf #1/\ }}

\catcode`@=12


\def\limpl{\mathrel{\Rightarrow}}
\def\repl{\mathrel{\Leftarrow}}

\def\land{\wedge}
\def\logand{\,\&\,}

\def\thinks{\models}


\def\pre#1#2{{}^{#1}\!#2}
\def\erp#1^#2{\pre{#2}{#1}}
\def\to{\rightarrow}

\def\cut{\cap}
\def\extend{\mathord{{}^{\frown}}}
\def\vxtend{\!\!\!\extend}

\def\Forces_#1``#2''{\forces_{#1}\hbox{``#2''}}

\def\fct{\pre{\omega}{\omega}}

\def\finSeq{\omega^{\lless\omega}}
\def\finxSeq#1{{#1}^{\lless \omega}}

\def\zerooneseq{\pre{\omega}{2}}
\def\finzerooneseq {2^{\lless \omega}}
\def\kzerooneseq#1 {\pre{#1}{2}}

\def\fsq#1^#2{\pre{#2}{#1}}  
\def\oob/{$\fct$-bounding}


\def\sqo#1:#2{\<{#1_{#2}:#2 < \omega _{1}}>}
\def\<#1>{\hbox{$\langle #1\rangle$}}  

\def\sqn#1{\<#1_n:n<\omega>}
\def\iter{\hbox{$\<P_\alpha, Q_\alpha: \alpha<\varepsilon >$}}

 \def\stronger{\ge}
 \def\weaker{\le}
 \def\emptycondition{\emptyset}
\def\lessdot{\mathrel{\mathord{<}\!\!%
  \raise 0.8 pt\hbox{$\scriptstyle\circ$}}}


\def\dom{\mathop{\rm dom}\nolimits}
\def\rng{\mathop{\rm rng}\nolimits}
\def\min{\mathop{\rm min}}

\def\sll{\rm}
\def\stem{{\sll stem}}
        \def\supp{\dom}
\def\succ{{\sll succ}}

\def\split{{\sll split}}
\def\height{{\sll ht}}
\def\sll{\rm}  

\def\vu{\nu}
\def\trim#1 #2{{#1}^{[#2]}}

\def\lless{\mathord{<}}
\def\lleq{\mathord{\le}}

\def\bool#1 {[\![ #1 ]\!]}

\long\def\set(#1:#2&#3){\setbox0=\hbox{#2}%
\{#1 :
\vtop{\hsize=\wd0\parindent=0cm\parskip=0cm%
\rightskip=0pt plus 0.1\wd0 minus 0.1\wd0%
#2 #3$
\}$}}




%


\baselineskip=1.08\baselineskip
\lineskiplimit=.92\lineskiplimit

\def\name#1{\mathchoice%
{\setbox0=\hbox{$\displaystyle #1$}
\setbox1=\vtop{\ialign{##\crcr
$\hfil{\displaystyle #1}\hfil$\crcr\noalign{\kern2pt%
\nointerlineskip}$\hfil\mathord{\displaystyle \sim}%
\hfil$\crcr\noalign{\kern3pt\nointerlineskip}}}%
\setbox2=\hbox{$\displaystyle \sim$}%
\box1}%
{\setbox0=\hbox{$\textstyle #1$}
\setbox1=\vtop{\ialign{##\crcr
$\hfil{\textstyle #1}\hfil$\crcr\noalign{\kern1.2pt%
\nointerlineskip}$\hfil\mathord{\textstyle \sim}%
\hfil$\crcr\noalign{\kern3pt\nointerlineskip}}}%
\setbox2=\hbox{$\textstyle \sim$}%
\wd1=\wd0\dp1=0cm\ifdim\wd2>\wd1 \wd1=\wd2\else\relax\fi
\ht1=\ht0\relax
\box1}{\setbox0=\hbox{$\scriptstyle #1$}
\vtop{\ialign{##\crcr
$\hfil{\scriptstyle #1}\hfil$\crcr\noalign{\kern1.4pt%
\nointerlineskip}$\hfil\mathord{\scriptstyle \sim}%
\hfil$\crcr\noalign{\kern2.1pt\nointerlineskip}}}%
}{\setbox0=\hbox{$\scriptscriptstyle #1$}
\vtop{\ialign{##\crcr
$\hfil{\scriptscriptstyle #1}\hfil$\crcr\noalign{\kern1pt%
\nointerlineskip}$\hfil\mathord{\scriptscriptstyle \sim}%
\hfil$\crcr\noalign{\kern1.5pt\nointerlineskip}}}%
}}

{\catcode`\@=11\gdef\neubox{\alloc@ 4\box \chardef \insc@unt}}
\def\assign#1{\neubox\current
\setbox\current=\hbox{$\name #1$}
\edef\next{\noexpand\edef
\csname n#1\endcsname{{\noexpand\copy\the\current}}}\next}

\def\nofirst#1{}
\def\assignp#1{\neubox\current
\setbox\current=\hbox{$\name #1$}
\edef\next{\noexpand\edef
\csname
n\expandafter\nofirst
\string#1\endcsname{{\noexpand\copy\the\current}}}\next}

\assignp\tau
\let\nt=\ntau
\assign f
\assign g
\assign q
\assign Q
\assign p
\assignp\sigma
\assignp\alpha
\assign r
\assign x
\assign y
\assign k


\def\compile#1#2{\neubox\current
\setbox\current=\hbox{$#2$}
\edef#1{\noexpand\copy\the\current}}

\def\ZFC/{\hbox{\rm ZFC}}


\font\fourteenbold=cmbx10 at 14.4 true pt

{
\fourteenbold

\parindent=0cm\parskip=0.7\baselineskip

\begingroup
\everypar={\hskip\parfillskip}    

\vglue 1 cm\relax

STRONG MEASURE ZERO SETS

WITHOUT COHEN REALS

\bigskip\bigskip\bigskip
{\rm January 1991}
\bigskip\bigskip
\parskip=\baselineskip

Martin Goldstern$^1$

{\rm Bar Ilan University}
\medskip

 Haim Judah$^1$

{\rm Bar Ilan University and  U.C.\ Chile}
\medskip

Saharon Shelah$^{1,2}$

{\rm Hebrew University of Jerusalem}

\endgroup
\bigskip\bigskip

\everypar={}

{\rm \leftskip=0.15\hsize\rightskip=0.15\hsize plus 1.5cm minus 1.5cm
\parindent=1true cm
\parskip=0cm
\noindent 
ABSTRACT.\quad  If ZFC is consistent, then each of the following are
                 consistent with ZFC + $2^{{\aleph_0}}=\aleph_2$: 
\endgraf
\advance\rightskip by20pt
\item{1.}  $X \subseteq \R$ is of strong measure zero iff $|X| \le
                 \aleph_1$ + there is a
                 generalized  Sierpinski set.
\item{2.} The union of $\aleph_1$ many strong measure zero sets is a
                 strong measure zero set +  there is a strong measure
                 zero set of size $\aleph_2$. 
\par
}

\vfill
\pageno=0\parskip=0cm
\rm $^1$  The authors thank the
Israel Foundation for Basic Research, Israel Academy of
Science.

$^2$ Publication 438

\footline={\hfil}
\eject

}

\advance\topskip by 10 pt
\headline={\headerfont
 \hfil
Goldstern, Judah, Shelah:  Strong measure zero sets
without Cohen reals \hfil
}
\font\headerfont=cmsl10 at 10 true pt
\footline={{\tt \jobname}\hfill{\tenrm \folio}\hfill{\tenrm\heute}}
\def\heute{January 1991}

\proofmodefalse
\ifproofmode\advance\hoffset by 0.7 true cm
\else 
\fi

\neuchap

{\bf \S0. Introduction}

In this paper 
   we continue the study of the structure of strong measure
   zero sets.  Strong measure zero sets have been studied from the
   beginning of this century.  They were discovered by E.~Borel, and
   Luzin, Sierpinski, Rothberger and others turned their attention to
   the structure of these sets and proved very interesting
   mathematical theorems about them.
 Most of the constructions of strong measure zero sets involve Luzin
   sets, which have a strong connection with Cohen reals (see
   [\smzRapid]).  In this  paper
    we will show that this connection is
   only apparent; namely, we will build models where there are strong
   measure zero sets of size $\c$ without adding Cohen reals over the
   ground model.

Throughout this work we will investigate questions about strong
   measure zero sets under the assumption  that
   $\c=2^{\aleph_0}=\aleph_2$.  The reason is that CH makes many of
   the    questions we investigate trivial, and there is no good
   technology available to deal with most of our problems when
   $2^{\aleph_0} > \aleph_2$.

\bigskip

\ppro\smzDef/Definition: A set $X \subseteq \R$ of reals has strong
   measure zero if for every sequence 
   $\<{\varepsilon}_i:i<\omega>$ of positive real numbers there is a
   sequence  $\<x_i: i<\omega>$ of real numbers such
   that 
    $$ X \subseteq \bigcup_{i<\omega} (x_i-{\varepsilon}_i,
                           x_i+{\varepsilon}_i)$$

   We let $\smz \subset \P(\R)$ be the ideal of strong measure zero
   sets.

\ppro\smzDefRem/Remark: (a) if  we work in $\zerooneseq$ 
   then $X \subseteq \zerooneseq$  has strong
   measure zero if  
    $$ (\forall h\in \fct )
           (\exists g \in \prod_n \kzerooneseq {h(n)} )
              (\forall x \in X)
                 (\exists^\infty n) (g(n)=x\on h(n))$$
  or equivalently, 
    $$ (\forall h\in \fct )
           (\exists g \in \prod_n \kzerooneseq {h(n)} )
              (\forall x \in X)
                 (\exists n) (g(n)=x\on h(n))\leqno (*)$$

(b) To  every question about strong measure zero sets in $\R$ there is
   a corresponding question about a strong measure zero set of
   $\zerooneseq$, and for all the questions we consider the corresponding
   answers are the same.    So we will work sometimes in $\R$,
   sometimes in $\zerooneseq$.

\ppro\defdef/Definition:   Assume that $\H \subseteq \fct $. We  say
that $\bar \nu  $ 
   has index $\H$, if $\bar \nu =  \<\nu^h:h\in \H>$ and  for all
   $h\in\H$, $\nu^h$ is a function with domain 
   $\omega$ and $\forall n\, \nu^h(n)\in \kzerooneseq h(n) $.  We let
   $$X_{\bar \nu}:= \bigcap_{h\in\H} \bigcup_{k\in \omega}
[\nu^h(k)]$$ 
 (where we let $[\eta]:= \{f\in \zerooneseq: \eta \subseteq
f\}$).

   We say that $X_{\bar \nu} $ is the set ``defined'' by $\bar \nu$.  

\ppro\defdefFact/Fact:  Assume $\H \subseteq \pre{\omega}{\omega}$  is a 
   dominating  family, i.e., for all $f\in\fct $ there is $h\in\H$
   such that $\forall n$ $f(n) < h(n)$.    Then:
\begindent
\ite 1  If $\bar \nu$ has index $\H$, 
   then $X_{\bar\nu}$ is a strong measure zero set.
\ite 2  If $X$ is a strong measure zero set, then there is a sequence
   $\bar\nu$ with index $\H$ such that  $X \subseteq X_{\bar\nu} $. 
\endent

\ppro\GluzinDef/Definition: A set of reals $X \subseteq \R$ is a
   GLuzin (generalized Luzin) set if for every meager set $M \subseteq
   \R$, $X \cut M$ has cardinality less than $\c$.  $X$ is a
   generalized Sierpinski set if set if for every set $M \subseteq \R$
   of Lebesgue measure 0, $X \cut M$ has cardinality less than $\c$.

\ppro\GluzinFact/Fact: (a) If $\c$ is regular, and $X$ is GLuzin, then
   $X$ has strong measure zero.

(b) A set of mutually independent Cohen reals over a model $M$ is a
   GLuzin set. 

(c) If $\c>{\aleph_1}$ is regular, and $X$ is a GLuzin set, then 
  $X$ contains Cohen reals over $L$. 

Proof: See [\smzRapid]. 

\ppro\onetheo/Theorem: [\smzRapid]  Con(ZF) implies Con(ZFC + there is a
   GLuzin set which is not strong measure zero). 

\ppro\twotheo/Theorem: [\smzRapid]  Con(ZF) implies Con(ZFC +
   $\c>{\aleph_1}$ + $\exists X \in [\R]^{\c}$, $X$ a strong measure
   zero set + there are no GLuzin sets). 

In theorem \bigsmz/ we will show a stronger form of \twotheo/. 

\ppro\SMZdef/Definition: 
\begindent
\ite 1 
Let \SMZ/ be the following statement:  ``Every
   set of reals of size less than $\c$ is a strong measure zero set.''
\ite 2  We say that the ideal of strong measure zero sets is
   $\c$-additive, or $\Add(\smz)$, if  for every $\kappa<\c$ the union
   of $\kappa$ many strong measure zero sets is a strong measure zero
   set.   (So $\Add(\smz) \limpl \SMZ/$.)
\endent

\ppro\RothsEq/Remark: 
Rothberger ([\EigenschaftC] and [\ProprieteC]) proved that the
   following are equivalent: 
  \begindent
\ite i \SMZ/
  \ite ii for every $h: \omega \to \omega $, for every $F\in [\prod_n
   h(n)]^{\lless \c}$, there exists $f^*\in \fct$ such that
   for every $f\in F$ there are infinitely many $n$ satisfying
   $f(n)=f^*(n)$. 
\endent

Miller ([\SomeProperties]) noted  that this implies the following:
   $$\Add(\meager)\quad\hbox{ iff }\quad\SMZ/\hbox{ and } \b=\c$$

(See \sNotation/   for  definitions)

\bigskip\bigskip

Rothberger proved interesting results about the existence of strong
   measure zero sets, namely: 

If $\b= \aleph_1 $, then there is a strong measure zero set of size
   $\aleph_1$.
(See [\JSW].)

In this spirit, we first prove the following result: 

\ppro\rothTheo/Theorem: If \SMZ/ and $\d=\c$, then there exists a
   strong measure zero set of size $\c$. 

 We start the proof by proving the following
\ppro\rothFact/Fact:  If $\d=\c$, then there is a set $\{ f_i:i<\c\}$
   of functions in $\fct$ such that for every $g\in\fct $, the set
   $$\{ i<\c: f_i\le^* g\}$$ 
   has cardinality  less than $\c$. 

Proof of the fact: We build $\<f_i:i<\c> $ by transfinite induction.
   Let $\fct = \{g_j:j<\c\}$.  We will ensure that for $j < i$,
    $f_i \not<^*g_j$.  This will be sufficient. 

    But this is easy to achieve, as for any $i$, the family
   $\{g_j:j<i\}$ is not dominating, so there exists a function $f_i$
   such that for all $j<i$, for infinitely many $n$, $f_i(n)>g_j(n)$. 

 This completes the proof of \rothFact/.

\ppro\rothTproof/Proof of \rothTheo/:  Using $\d=\c$, let $\<f_i:i<\c>$
   be a sequence as in \rothFact/. Let $F:    \fct \to [0,1]-\Q$ be a
   homeomorphism. ($\Q$ is the set of rational numbers.) 
   We will show that $X:=    \{F(f_i): i < \c\} $ is a
   strong measure zero set.  

   Let $\sqn \varepsilon$ be a sequence of positive numbers.  Let
   $\{r_n:n\in \omega \}=\Q$.   Then $U_1:= \bigcup_{n\in \omega }
   (r_n - \varepsilon_{2n}, r_n + \varepsilon_{2n})$
   is an open set. So  $K:= [0,1]-U_1$ is closed, hence compact.   As
   $K \subseteq \rng(F)$, also $F^{-1}(K) \subseteq \fct$ is a compact
   set.   So for all $n$ the projection of $F^{-1}(K)$ to the
   $n$th coordinate is a compact (hence bounded) subset of $\omega$,
   say $ \subseteq g(n)$. So 
      $$ F^{-1}K \subseteq \{f \in \fct: f \le^* g\}$$

    Let $Y:= X - U_1 \subseteq K $.  Then  $Y \subseteq F(F^{-1}(K))
   \subseteq \{F(f_i): f_i \le^* g\}$ is 
   (by assumption on $\<f_i:i<\c>$) a set of size $<\c$, hence has
   strong measure zero.  So there exists a sequence of real numbers
   $\sqn x $ such that $Y \subseteq U_2$, where $$ U_2:= 
    \bigcup_{n\in \omega }
   (x_n - \varepsilon_{2n+1}, x_n + \varepsilon_{2n+1})$$ and $X
   \subseteq U_1 \cup U_2$.   So $X$ is indeed a strong measure zero
   set. 
\medskip

\bigskip

In section 2 we will build models where $\Add(\smz)$ holds and the
   continuum is bigger than $\aleph_1$ without adding Cohen reals.
First we will show in \smallproof/: 

\medskip

\ppro\smallsmz/Theorem:  If \ZFC/ is consistent, then 

   \centerline{ \ZFC/ + $\c=\aleph_2$ + $\smz=[\R]^{\lleq\aleph_1 }$ +
   there are no Cohen reals over $L$}
   is consistent. 

Note that $\c=\aleph_2$ and $\smz=[\R]^{\lleq\aleph_1}$ implies 
 \begindent
\ite 1  $\Add(\smz)$.  (Trivially)
 \ite 2 $\b=\d=\aleph_1$. (By \rothTheo/)
\endent
      
The same result was previously obtained by Corazza[\Corazza].  In his
model the nonexistence of strong measure zero sets of size $\c$ is
shown by proving that every set of size $\c$ can be mapped
uniformly
continuously onto the unit interval (which is impossible for a strong 
measure zero set).  Thus, the question arises
whether is possible to get a model of $\smz=[\R]^{<\c}$ +
$\c=\aleph_2$ + ``not all set of size $\c$ can be 
continuously mapped onto $[0,1]$.''

By adding random reals to our construction, we answer this question in
the following stronger theorem:

\ppro\smallprime/Theorem:  If \ZFC/ is consistent, then 

   \centerline{ \ZFC/ + $\c=\aleph_2$ + $\smz=[\R]^{\lleq\aleph_1 }$ +
   there are no Cohen reals over $L$ }
\centerline{ +  there is a generalized Sierpinski set}
   is consistent. (See \GluzinDef/.)

By a remark of Miller [\Mapping, \S2] a generalized Sierpinski set
cannot be mapped continuously onto $[0,1]$ (not even with a Borel
function).

\bigskip

Pawlikowski [\Pawlikowski] showed that the additivity of the ideal
   $\smz$ of strong measure 
   zero sets does not imply the additivity of the ideal $\meager$ of
   meager sets.  For this he built a model satisfying 
   $\Add(\smz)$ + $\c=\aleph_2 $ + $ \b=\aleph_1$. He used a finite
   support iteration    of length $\omega_2$.  So he adds many Cohen
   reals, and in the final model $\Cov(\meager) $ holds (i.e., $\R$
   can not be written as the union of less than $\c$ many meager
   sets).   We will improve his result in the next theorem:

\ppro\bigsmz/Theorem:  If \ZFC/ is consistent, then 

   \centerline{\ZFC/ + $\c=\d=\aleph_2>\b$ + $\Add(\smz)$ 
       +  no real is Cohen over $L$}
   is consistent. 

(Note that by \rothTheo/, $\d=\c$ + $\Add(\smz)$ implies that there
   is a strong measure zero set of size $\c$.)

\ppro\sNotation/Notation:  We use standard set-theoretical notation. 
   We identify natural
   numbers $n$ with their set of predecessors,  $n=\{0, \ldots, n-1\}$.
   $\pre AB$ is the set of functions from $A$ into $B$,
   $\finxSeq{A}:= \bigcup_{n<\omega} \pre n A $.
   $|A|$ denotes the cardinality of a set $A$.   $\P(A)$ is the power
   set of a set $A$, $A \subset B$ means $A \subseteq B\logand
   A\not=B$.  $A-B$ is the set-theoretic difference of $A$ and $B$. 
   $[A]^{\kappa}:= \{X \subseteq A: |X|=\kappa\}$.  $[A]^{\lless
   \kappa}$ and $[A]^{\lleq \kappa}$ are defined similarly. 
   (We write $A:=B $ or $B=:A$ to mean:  the expression $B$ defines the
   term (or constant) $A$.)

   $Ord$ is the set of ordinals. $cf(\alpha)$ is the cofinality of an
   ordinal $\alpha$.  $S^\alpha_{\beta} :=\{{\delta}\in
   \omega_{\beta}: cf({\delta})=\omega_\alpha\}$.  In particular,
   $S^1_2$ is the set of all ordinals $<\omega_2$ of uncountable
   cofinality.

   $\R$ is the set of real numbers.    $\c = |\R|$ is the size of the
   continuum.   For $f,g\in \fct$ we let $f<g$ iff for all $n$
   $f(n)<g(n)$,  and $f<^* g $ if for some $n_0\in 
   \omega$, $\forall n\ge n_0$ $f(n)<g(n)$.   The ``bounding number''
   $\b$ and the ``dominating number'' $\d$ are defined as 
   $$\eqalign{\b := &\min\{|\H|: \H \subseteq \fct, \forall g\in \fct\,
     \exists h\in \H \,\,\lnot(h<^*g)\}\cr
             \d:= &\min\{|\H|: \H \subseteq \fct, \forall g\in \fct\,
     \exists h\in \H \,\, g<h\}\cr
       = & \min\{|\H|: \H \subseteq \fct, \forall g\in \fct\,
     \exists h\in \H \,\, g<^*h\}\cr }$$
   (It is easy to see $\omega_1 \le \b \le \d \le \c$.)

We call a set $\H \subseteq \fct$ dominating, if 
  $\forall g\in \fct  \,    \exists h\in \H \,\, g<h$.

$\meager$ is the ideal of meager subsets of  $\R$  (or of
   $\zerooneseq$). 
$\smz$ is the ideal of strong measure zero sets.     For any ideal $\J
   \subset \P(\R)$,  $\Add(\J)$ abbreviates the statement: 
   ``The union of less than $\c$ many sets in $\J$ is in $\J$.'' 
   $\Cov(\J)$ means that the reals cannot be covered by less than $\c$
   many sets in $\J$.

If $f$ is a function, $\dom(f)$ is the domain of $f$, and $\rng(f)$ is
the range of $f$. For   $A \subseteq \dom(f)$, $f\on A$ is the
restriction of $f$ to $A$.   For $\eta\in \finzerooneseq$,
$[\eta]:=\{f\in \zerooneseq: \eta \subseteq f\}$.

\ppro\mNotation/More Notation:
If $Q$ is a forcing notion, $G_Q$ is the canonical name for the
   generic filter on $Q$.   We interpret $p \weaker q$  as $q$ is
   stronger (or ``has more information'') than $p$.  (So $p \weaker q
   \limpl q \forces p\in G_Q$.)

When we deal with a (countable support) iteration $\<P_\alpha,
   Q_\alpha:\alpha<\varepsilon>$, we write $G_\alpha$ for the
   canonical name of the generic filter on $\alpha$, and $G(\alpha) $
   for the generic filter on $Q_\alpha$.
If there is a natural way to associate a ``generic'' real to the
generic filter on $Q_\alpha$, we write $g_\alpha$ for the real given
by $G(\alpha)$.   We write
   $\forces_\alpha$ for the forcing relation of $P_\alpha$.   
 If
   ${\beta}<\alpha$, $G_{\beta}$ always stands for $G_\alpha\cut
   P_{\beta}$.  $V= V_0$ is the ground model, and $V_\alpha =
   V[G_\alpha]$.  $P_\varepsilon$ is the countable support limit of
   $\<P_\alpha:\alpha<\varepsilon>$.  $P_\varepsilon/G_\alpha$ is the
   $P_\alpha$-name for $\{p\in P_\varepsilon: p\on \alpha\in
   G_\alpha\}$ (with the same $\le$ relation as $P_\varepsilon$).  The
   forcing relation with respect to $P_\varepsilon/G_\alpha$ (in
   $V_\alpha$) is denoted by $\forces_{\alpha \varepsilon}$.

There is a natural dense embedding from $P_\varepsilon $ into
$P_\alpha * P_\varepsilon/G_\alpha$.   Thus  we  always identify
$P_\alpha$-names for $P_\varepsilon/G_\alpha$-names with the
corresponding $P_\varepsilon$-names.

  $\emptycondition_\alpha$
is the weakest condition of $P_\alpha$, and $\emptycondition_\alpha
\forces_\alpha \varphi $ is usually abbreviated to $\forces_\alpha
\varphi$.   (So $\forces_\alpha(\forces_{\alpha {\delta}} \varphi)$
iff $\forces_{\delta} \varphi$).

\bigskip

\ppro\emNot/Even more Notation:
\begingroup
\parskip=0cm
The following notation is used when we deal with trees of finite
   sequences: 
\smallskip
 For $\eta\in V^{\lless \omega}$, $i\in V $, $\eta\extend
   i$ is the 
     function $\eta\cup \{\<|\eta|, i>\}\in V^{\lless \omega}$. 

 $p \subseteq \finSeq $ is a tree if $p\not=\emptyset$, and for
     all $\eta\in p$, all $k<|\eta|$, $\eta\on k \in p$.   Elements of
a tree are often called ``nodes''.   We call $|\eta|$ the ``{length}''
   of $\eta$. {\bf We reserve the word  ``{height}'' for the notion
   defined in  \splitDef/}.

For $p \subseteq \finSeq$, $\eta\in p$, we let 
    $ \succ_p(\eta):=\{i: \eta\extend i \in p\}$.   

 If $p$ is a tree,  $\eta\in p$, let
 $ \trim p \eta := \{\nu\in p : \eta \subseteq \nu \hbox{ or }
				  \nu \subseteq \eta \} $. 

 If $p \subseteq \finSeq$ is a tree,  $b \subseteq p $ is
          called a  {\it branch}, if $b$ is 
          a maximal subset of $p$ that is linearly ordered by
          $\subseteq$.  

Clearly, if $\forall \eta\in p\,\,\succ_p(\eta)\not=\emptyset$, then a
   subset $b \subseteq p$ is a branch iff $b$ is of the form $b=\{f\on
   n:n\in \omega\}$ for some $f\in \fct$.

We let $\stem(p)$
   be the intersection of all branches of $p$. 
\endgroup

\bigskip

\par


\bigskip\goodbreak\bigskip

{\bf \S1.   A few well known facts}
\neusection

We collect a few more or less well known facts about forcing, for
later reference. 

\ppro\ppDef/Definition: An ultrafilter $\U$ on $\omega$ is called a
{\bf P-Point}, if for any sequence  $\<A_n:n\in\omega>$ of sets in
$\U$ there is a set $A$ in $\U$ that is almost contained in every
$A_n$ (i.e., $\forall n\,\, A-A_n$ is finite).

\ppro\ppGame/Definition:  For any ultrafilter $\U$ on $\omega$, we
   define the P-point game $G(\U)$ as follows: 

\begingroup\advance\leftskip by 2 cm\advance\rightskip by 2 cm
There are two
   players, ``IN'' and ``NOTIN''.   The game consists of  $\omega$ many
   moves. 
  
   In the $n$-th move,
   player NOTIN picks a set $A_n\in \U$, and player IN picks a finite set
   $a_n \subseteq A_n$. 

  Player IN wins if after $\omega $ many moves, $\bigcup_n a_n\in \U$.

  We write a play (or run) of $G(\U)$ as $$\<A_0; a_0 \to A_1; a_1\to
   A_2; \ldots\,>.$$

\endgroup

It is well known   that an ultrafilter
$\U$ is a P-point iff player 
NOTIN does not   have   a winning  strategy in the P-point game.

\begingroup
\def\\#1,#2.{[n_{#1}, n_{#2})} 
For the sake of completeness, we give a proof of the nontrivial
   implication ``$\limpl$'' (which is all we will need later):

Let $\U$ be a P-point, and let ${\sigma}$ be a strategy for player
   NOTIN. We will construct a run of the game in which player NOTIN
   followed ${\sigma}$, but IN won. 

Let $A_0$ be the first move according to ${\sigma}$. For each $n$, let
   $\A_n$ be the set of all responses of player notin according to
   ${\sigma}$ in an initial segment of a play of length $\le n$ in
   which player IN has played only subsets of $n$: 
  
$$ \A_n:=\set(A_k:
   $k\le n$, $\<A_0; a_0\to A_1; \ldots; a_{k-1}\to A_k>$ 
   is an &initial segment of a play in which NOTIN obeyed ${\sigma}$,
   and $a_0, \ldots, a_{k-1} \subseteq n$)$$

Note that $\A_0=\{A_0\}$, and for all $n$,  $\A_n$ is  a finite subset
of $\U$. 

As $\U$ is a P-point, there is a set $X\in\U$ such that for all $A\in
   \bigcup_n \A_n$, $X-A$ is finite.  

Let $X \subseteq A_0 \cup n_0$, and for $k>0$ let $n_k$ satisfy 
 $$ n_k > n_{k-1} \hbox{\quad and\quad} \forall A\in \A_{n_{k-1}}\,\,
   X \subseteq A \cup n_k $$
Either $\bigcup_{k\in \omega} \\2k,2k+1. \in \U$, or $\bigcup_{k\in
   \omega} \\2k+1,2k+2. \in \U$.

  Without loss of generality we assume
   $\bigcup_{k\in \omega} \\2k,2k+1.\in \U$.

Now define a play $\<A_0; a_0 \to A_1; a_1\to A_2; \ldots\,>$ of the
   game $G(\U)$ by  induction    as follows: 
\begindent
\item{}  $A_0$ is given.
\item{} Given $A_j$, let $a_j:= A_j\cut \\2j,2j+1.$ and
   let $A_{j+1}$ be ${\sigma}$'s response to $a_j$. 
\endent
Then as $a_0$, \dots, $a_{j-1} \subseteq n_{2j}$, we have $X \subseteq
   A_j \cup n_{2j}$ for all $j$.  Therefore for all $j$ we have
   $X\cut \\2j,2j+1.  \subseteq  (A_j\cup n_{2j} )\cut \\2j,2j+1.  =
   A_j \cut \\2j,2j+1. = a_j$.  So $\bigcup _{j\in \omega}a_j \supseteq X
   \cut \bigcup_{j\in \omega} \\2j,2j+1. \in \U$.

Thus  player IN wins
   the play $\<A_0; a_0\to A_1; a_1\to A_2; \ldots >$ in which player
   NOTIN obeyed ${\sigma}$. 
\endgroup

\bigskip

\ppro\pppDef/Definition:  
We say that a forcing notion $Q$ preserves P-points, if for every
   P-point ultrafilter $\U$ on $\omega$, $\forces_Q$``$\U$ generates
   an ultrafilter'', i.e. $\forces_Q$ $\forall x\in \P(\omega)\,
   \exists u\in\U\,\,(u \subseteq x \hbox{ or } u \subseteq \omega-x)$.''

\medskip
[\Rational] defined the  following forcing notion: 

\ppro\RPSdef/Definition: ``Rational perfect set forcing'', $RP$ is
defined as the set of trees $p \subseteq \finSeq$ satisfying
\begindent
\ite 1        for all $\eta \in p$, $|\succ_p(\eta)|\in
   \{1,{\aleph_0}\}$ (See \emNot/)
\ite 2  for all $\eta\in p$ there is $\vu\in p$ with $\eta \subseteq
   \nu$ and $|\succ_p(\eta)| = {\aleph_0}  $.
\endent
We let $p \stronger q $ iff $p \subseteq q$.

Then the following hold: 

\ppro\MillersLemma/Lemma: 
\begindent
\ite 1  $RP$ preserves P-points. ([\Rational, 4.1])
\ite 2 $RP$ adds an unbounded function. ([\Rational, \S2])
\ite 3 $RP$ is proper. (This is implicit in [\Rational]. See also
   \RPproper/) 
\endent

 \medskip
The next lemma can be found, e.g., in [\Kunen, VII ?? and Exercise H2]: 

\ppro\ccFact/Fact:  If $Q$ is a forcing notion satisfying the
   $\aleph_2$-cc, then 
\begindent
\ite  1 If $\forces_Q \name c:\omega_2^V\to \omega_2^V$, then there is a
   function $c:\omega_2\to {\omega_2} $ such that $\forces_Q \forall
   \alpha<\omega_2: \name c(\alpha)<c(\alpha)$. 
\ite 2 $\forces_Q \aleph_2^V = \aleph_2$. 
\ite 3 For every stationary $S \subseteq \aleph_2$,
   $\forces_Q$ ``$S$ is stationary on~$\aleph_2$''. 
\endent

\medskip

The following fact is  from [\ProperForcing, V 4.4]: 
\ppro\nonew/Fact: 
Assume $\<P_\alpha, Q_\alpha:\alpha<\omega_2> $ is an iteration of
proper forcing notions $Q_\alpha$. Then for every  ${\delta}\le
\omega_2$ of cofinality $>\omega$, $\forces_{\delta}\, \fct \cut
V_{\delta} = \fct \cut \bigcup_{\alpha<{\delta}} V_\alpha$, or in
other words: ``no new reals appear in limit stages of cofinality
$>\omega$''.  

As a consequence, $\forces_{\omega_2}$``If $X \subseteq \fct$, $|X|\le
\aleph_1$, then there is ${\delta}<{\omega_2} $ such that $X\in
V_{\delta}$.''

\medskip

We also recall the following facts about iteration of proper forcing
              notions: 

\ppro\ccLemma/Lemma: Assume $CH$, and let $\<P_\alpha,
   Q_\alpha:\alpha<\omega_2>$ be a countable support iteration such
   that for all $\alpha<\omega_2$, $\forces_\alpha$``$Q_\alpha$ is a
   proper forcing notion of size $\le \c$.''   

Then
\begindent
\ite 1 $\forall \alpha<\omega_2$: $\forces_\alpha \c=\aleph_1$.
              (see [\ProperForcing, III 4.1]) 
\ite 2 $\forces_{\omega_2} \c \le \aleph_2$. (This follows from
              \nonew/ and (1))
\ite 3 For all $\alpha\le {\omega_2} $, $P_\alpha$ is proper
              [\ProperForcing, III 3.2] and  satisfies the
              $\aleph_2$-cc. (See [\ProperForcing, III 4.1])
\ite 4 $\forces_{{\omega_2}}\aleph_1^V=\aleph_1$. (See
   [\ProperForcing, III 1.6]) 
\endent
\medskip

In [\BlassShelah, 4.1]  the following is proved:

\ppro\preservePP/Lemma: Assume $\<P_\alpha, Q_\alpha:\alpha<\omega_2>$
   is as in \ccLemma/, and for all $\alpha<\omega_2$: 
    $$ \forces_\alpha \hbox{``$Q_\alpha$ preserves P-points.''}$$
 Then for all $\alpha\le {\omega_2}$, $P_\alpha$ preserves P-points. 

\medskip

\ppro\oobDef/Definition:  We say that a forcing notion $Q$ is \oob/,
   if the set of ``old'' functions is a dominating family in the
   generic extension by $Q$, or equivalently, 

   $$ \forces_Q \forall f\in \fct\,\exists g\in \fct\cut V \,\,
   \forall n\, f(n)<g(n)$$

 [\ProperForcing, V 4.3] proves:

\ppro\preserveoo/Lemma: Assume $\<P_\alpha, Q_\alpha:\alpha<\omega_2>$
   is as in \ccLemma/, and for all $\alpha<\omega_2$: 
    $$ \forces_\alpha \hbox{``$Q_\alpha$ is \oob/ and $\omega$-proper.''}$$
 Then for all $\alpha\le{\omega_2}$, $P_\alpha$ is \oob/.

(We may even replace $\omega$-proper by ``proper'', see
[\ProperForcing], [\ShelahsPreservation])

The following is trivial to check: 

\ppro\noCfact/Fact:  Assume $Q$ is a forcing notion that preserves
   P-points or is \oob/.  Then 
    $$\forces_Q \hbox{``There are no Cohen reals over $V$''}$$

\medskip

\ppro\strongoo/Definition:  A forcing notion $P$ is {\bf strongly
\oob/}, if there is a sequence $\<{\le}_n:n\in \omega>$ of binary
reflexive    relations on $P$ such that for all $n\in \omega$:
\begindent
\ite 1 $p \le_n q \,\, \limpl \,\, p \le  q$. 
\ite 2 $p \le_{n+1} q \,\, \limpl \,\, p \le_n q$. 
\ite 3 If $p_0 \weaker_0 p_1 \weaker_1 p_2 \weaker_3 \cdots\,$, then
there is a $q$ such that $\forall n \,\,  p_{n+1} \weaker_n q $. 
\ite 4  If $p\, \forces $``$\nalpha$ is an ordinal,'' and $n\in \omega$,
        then there exists $q \stronger_n p$ and a finite set $A
        \subseteq Ord$ such that $Q \forces \nalpha\in A$. 
\endent

\ppro\strongooIt/Definition: (1) If \iter{} is an iteration of strongly
        \oob/ forcing notions, $F \subseteq \varepsilon $ finite,
        $n\in \omega$, $p,q\in P_\varepsilon$, we say that $p
        \weaker_{F,n} q$ iff $p \weaker q$ and $\forall \alpha\in F\,\,
        q\on \alpha \forces p(\alpha) \weaker_n q(\alpha)$. 

(2) A sequence $\<{\<p_n, F_n>}: n\in \omega>$ is called a fusion
        sequence if $\<F_n:n\in \omega>$ is an increasing family of
        finite subsets of $\varepsilon$, $\<p_n:n\in \omega>$ is an
        increasing family of conditions in $P_\varepsilon$, $\forall
        n\,\, p_n \weaker_{n,F_n} p_{n+1}$ and $\bigcup_n \dom(p_n)
        \subseteq \bigcup_n F_n$.

Note that \strongoo/ is not a literally a strengthening of
        Baumgarter's 
        ``Axiom A'' (see [\IteratedForcing]), as we do not require
        that the relations $\le_n$ are transitive, and in (2) we only
        require $p_{n+1} \weaker_n q$ rather than $p_{n+1}
        \weaker_{n+1} q$. Nevertheless, the same proof as in
        [\IteratedForcing] shows the following fact:

\ppro\strooFact/Fact: 
\begindent
\ite 1 If the  sequence $\<{\<p_n, F_n>}: n\in \omega>$ is  a fusion
        sequence, then there exists a condition $q\in P_\varepsilon$
        such that for all $n\in \omega$, $p_{n+1} \stronger_{F_n,n}
        q$. 
\ite 2 If $\nalpha$ is a $P_\varepsilon$-name of an ordinal, $n\in
        \omega$, $F \subseteq P_\varepsilon$ finite, then for all $p$
        there exists a condition $q \stronger_{n,F} p$ and a finite
        set $A$ of ordinals such that $q \forces \nalpha \in A$. 
\ite 3 If $\name X$ is a $P_\varepsilon$-name of a countable set of
           ordinals, $n\in \omega$, $F \subseteq
        P_\varepsilon$ finite, then for all $p$ there exists a 
        condition $q \stronger_{n,F} p$ and a countable set $A$ of
        ordinals such that $q \forces \name X \subseteq  A$.
\endent

The next fact is also well known: 

\ppro\Random/Fact: Let $B$ be the random real forcing.  Then $B$ is
        strongly \oob/. 

[Proof: Conditions in $B$ are Borel subsets of $[0,1]$ of positive
        measure, $p \weaker q $ iff $ p \supseteq q$.    We let $p
        \weaker_n q$ iff $p \weaker q $ and $\mu(p-q) \le
        10^{-n-1}\mu(p)$, where $\mu$ is the Lebesgue measure.  Then if
        $p_0 \stronger_0 p_1 \stronger_1 \cdots$, letting
        $q:=\bigcap_n p_n$ we have for all $n$, all $k\ge n$,
        $\mu(p_k-p_{k+1} )\le 
        10^{-k-1}\mu(p_k )\le 10^{-k-1}\mu(p_n)$, so $\mu(p_n-q)\le
        10^{-n-1} + 10^{-n-2} + \cdots \,\le 2*10^{-n-1}\mu(p_n)$.  In
        particular, $\mu(q)\ge 0.8*\mu(p_0)$, so $q$ is a condition,
        and $q \stronger_{n-1} p_n$ for all $n>0$. 

      Given a name $\nalpha$, an integer $n$  and a condition $p$ such
        that $p \forces $``$\nalpha$ is  an ordinal,'' let $A$ be the
        set of all ordinals 
        ${\beta}$         such that $\bool \nalpha={\beta} \cut p$ has
        positive measure ($\bool \varphi {} $ is the boolean value of
        the statement $\varphi$, i.e. the union of all conditions
        forcing $\varphi$).  Since $\sum_{{\beta}\in A} \mu(\bool
        \nalpha={\beta} \cut p) = \mu(p)$ there is a finite subset $F
        \subseteq A$ such that letting $q:=p\cut \bigcup_{\beta\in A}
        \bool  \nalpha={\beta} $ we have $\mu(q)\ge
        (1-10^{-n-1})\mu(p)$.  So $q \stronger_n p$ and $q \forces
        \nalpha\in F$.]

\bigskip\bigskip

We will also  need the following lemma from [\RealValued, \S5, Theorem 9]: 

\ppro\splitstat/Lemma:  Every stationary $S \subseteq \aleph_2$ can be
   written as a union of $\aleph_2$ many disjoint stationary sets. 

\bigskip
\bigskip

Finally, we will need the following easy fact (which is  true for any
   forcing notion $Q$)

\ppro\rzeroone/Fact:   If $\nf$ is a $Q$-name for a function from
   $\omega$ to $\omega$, $\forces_Q\, \nf\notin V$, and $r_0$, $r_1$
   are any two 
   conditions in $Q$, then there are  $l\in\omega$,
   $j_0\not=j_1$, $r_0' \stronger r_0$, $r_1' \stronger r_1$ such that
   $r_0' \forces \nf(l)=j_0$, $r_1' \forces \nf(l)=j_1$. 

[Proof: There are a function $f_0$ and a sequence $r_0 =r^0 \weaker r^1
   \weaker \cdots$ of conditions in $Q$ such that for all $n$, $r^n
   \forces \nf\on n = f_0\on n$.  Since $r_1 \forces \nf \notin V$,
   $r_1 \forces \exists l\, \nf(l)\not=f_0(l)$. 
   There is a condition $r_1' \stronger r_1$ such that for some $l\in
   \omega$ and some $j_1\not=f_0(l)$, $r_1' \forces \nf(l)=j_1$.  Let
   $j_0:=f_0(l)$, and let $r_0' := r^{l+1}$.]

\bigskip\goodbreak\bigskip


\def\ivec{{\vec{\imath}}}

{\bf \S2\quad H-perfect trees}
\neusection

In this section we describe a forcing notion ${PT_H}$ that we will use
   in an 
   iteration in the next section.  We will prove  the following
   properties of ${PT_H}$: 
\begindent
\ite a  ${PT_H}$ is proper and  \oob/.
\ite b  ${PT_H}$ preserves P-points. 
\ite c ${PT_H}$ does not ``increase'' strong measure zero sets defined
   in the ground model. 
\ite d   ${PT_H}$ makes the reals of the ground model (and hence, by (c),
        the union of all strong measure zero sets defined in the ground
        model) a strong measure zero set.
\endent
      
\bigskip

  \ppro\HPerfDef/Definition:  For each function $H$ with domain
   $\omega$ satisfying $\forall n\in \omega$ $1<|H(n)|<\omega$,  we
   define the forcing ${PT_H}$, the set of $H$-perfect trees 
     to be the set of all $p$ satisfying
  \begindent
  \ite A $ p \subseteq \finSeq $ is a tree. 
  \ite B $\forall\eta\in p\,\,\forall l\in\dom(\eta): 
     \eta(l) \in H(l)$. 
  \ite C $\forall \eta \in p: |\succ_p(\eta)| \in 
	\{1,\left|H\left(\left|\eta\right|\right)\right|\}$. 
  \ite D $\forall \eta\in p\, \exists\nu \in p: \eta \subseteq \nu, 
    | \succ_p(\nu)|= |H(|\nu|)|$. 
  \endent

\medbreak

  \ppro\splitDef/Definition:
  \begindent
  \ite 1   For $p\in {PT_H}$, we let the set of
     ``splitting nodes'' of $p$ be 
     $$ \split(p) := \{\eta\in p : |\succ_p(\eta)|>1\}$$
  \ite 2 The {\bf height} of a node $\eta\in p\in {PT_H}$ is  the number
of splitting nodes strictly  below $\eta$:
     $$\height_p(\eta):=|\{\nu\subset\eta:\nu\in\split(p)\}|$$ 
\item{}  (Note that $\height_p(\eta) \le |\eta|$.)
  \ite 3 For $p\in {PT_H}$, $k\in \omega$, we let the $k$th
          splitting level of $p$ be the set of splitting nodes of
        height  $k$. 
     $$ \split_k(p) := \{\eta\in \split (p) : \height_p(\eta)=k\}$$
\item{} (Note that $\split_0(p) =   \{\stem(p)\}$.) 
  \ite 4  For  $u \subseteq \omega$, we let 
        $$ \split^u(p) := \bigcup_{k\in u} \split_k(p)$$
  \endent

\ppro\fPerfRem/Remarks: 
\begindent
\ite i  Since  $H(n)$ is finite, (3) just means that either
   $\eta$ has a unique successor $\eta\extend i$, or 
     $\succ_p(\eta)= H(|\eta|)$.) 
\ite ii  Letting $H'(n)= |H(n)|$, clearly ${PT_H}$ is isomorphic to
   $PT_{H'}$ (and the obvious isomorphism respects the functions
   $\eta\mapsto \height_p(\eta)$, $\<p,k>\mapsto \split_k(p)$, etc)
\endent

\ppro\RPremark/Remark:  If we let $H(n)=\omega$ for all $n$, then 
\HPerfDef/(A)--(D) define $RP$, rational perfect set forcing. 
  The
definitions in \splitDef/ make sense also for this forcing.  Since we
will not use the fact that $H(n)$ is finite  before \finSplit/,
 \splitFac/--\joinFact/ will be true also for $RP$.

  \ppro\splitFac/Fact: Let $p,q\in {PT_H} $, $n\in \omega$, $\eta,\nu\in
\finSeq$.  Then 
\begindent
\ite 1 If $\eta \subset \nu\in p$, then $\height_p(\eta) \le
\height_p(\nu)$.   If moreover $\eta\in \split(p)$, then
$\height_p(\eta) < \height_p(\nu)$. 
\ite 2 If $b \subseteq p $ is a branch, then $b\cut
		\split_n(p)\not=\emptyset$. 
\ite 3 If $p \supseteq q$, then for all $n$, $q\cut
       \split_n(p)\not=\emptyset$. 
\ite 4 If $\eta\in p$ and $\height_p(\eta)\le n$ then $\exists \nu\in
       p$, $\eta \subseteq \nu$ and $\nu\in \split_n(p)$. 
\ite 5 If $\eta_0\not=\eta_1$ are elements of $\split_n(p)$, then
       $\eta_0\not\subseteq\eta_1$, and $\eta_1\not\subseteq\eta_0$.  
\endent

Proof: (1) is immediate form the definition of $\height$. 

For (2), it is enough to see that $b\cut \split(p)$ is infinite. (Then
ordering $b$ by inclusion, the $n$th  element of $b\cut \split(p)$
will be in $\split_{n-1} (p)$.) 

So assume that $b\cut \split(p)$ is finite.   Recall that each
$\eta\in b-\split(p)$ has a unique successor in $p$.   By
\HPerfDef/(C), $b$ cannot have a last element, so $b$ is infinite.
Hence there is ${\eta_0}\in b$ such that  
$$ \forall \vu\in b: {\eta_0} \subseteq \nu \limpl |\succ_p(\nu)|=1.$$
A trivial induction on $|\nu|$ shows that this implies
$$ \forall \nu \in p: {\eta_0} \subseteq \nu \limpl  \nu \in b. $$ 
Hence
$$ \forall \nu \in p: {\eta_0} \subseteq \nu \limpl  |\succ_p(\nu)|=1.$$
This contradicts \HPerfDef/(D).

To prove (3), let $b$ be any branch of $q$.  $b$ is also a branch of
       $p$, so (2) shows that $q\cut \split_n(p) \supseteq b \cut
       \split_n(p) \not=\emptyset$. 

Proof of (4): Let $b$ be a branch of $p$ containing $\eta$.   By (2)
       there is $\vu\in b\cut \split_n(p)$.  If $\nu \subset \eta$,
       then $\height_p(\eta)>\height_p(\nu)=n$, which is impossible.
       Hence $\eta \subseteq \nu$.

(5) follows easily from (1). 
\bigskip

  \ppro\leDef/Definition: For $p, q\in {PT_H}$, $n \in \omega$,   we let 
  \begindent
  \ite 1 $ p \weaker q $ (``$q$ is stronger than $p$'') iff $q \subseteq
     p$. 
  \ite 2 $p \weaker_n q $ iff $p \weaker q $ and $\split_n(p)
   \subseteq q $.  (So also $\split_k(p) \subseteq q $ for all $k<n$.)
  \endent

\ppro\sameStem/Fact:  If $p \weaker_n q$, $n>0$, then
        $\stem(p)=\stem(q)$.

\ppro\lenFact/Fact: Assume $p,q\in {PT_H}$, $n\in \omega$, $p
\weaker_n q$.  
\begindent
\ite 0 For all $\eta\in q$, $\height_q(\eta) \le \height_p(\eta)$. 
\ite 1 For all $k\le n$, $\split_k(p) \subseteq \split(q)$. 
\ite 2 For all $k<n$, $\split_k(p) = \split_k(q)$. 
\ite 3 If $p \weaker_n q \weaker_n r$, then $p \weaker_n r$. 
\endent

Proof: (0) is clear. 

(1): Let $\eta\in \split_k(p)$ for some $k<n$, then by
       \splitFac/(4) there is a $\nu $, $\eta \subseteq \vu \in
       \split_n(p) \subseteq q$, so $\eta\in q$. 

(2): Let $\eta\in \split_k(p)$, then $\eta\in \split(q)$. Clearly
       $\height_q(\eta) \le \height_p(\eta)=k$.  Using (1)
       inductively, we also get  $\height_q(\eta)\ge k$. 

(3): Let $\eta\in \split_n(p)$.  So $\eta\in q$, $\height_q(\eta) \le
       \height_p(\eta) = n$. By \splitFac/(4), there is $\nu\in
       \split_n(q)$, $\eta \subseteq \nu$.  As $\nu\in r$, $\eta\in
       r$.

\ppro\fusionFact/Definition and Fact: If 
$p_0 \weaker_1 p_1 \weaker_2 p_2 \weaker_3 \cdots\,$ are conditions in
     ${PT_H}$, then we call the sequence $\sqn p $ a ``fusion
sequence''.   If 
          $\sqn p $ is a fusion sequence, then 
  \begindent
  \ite 1 $p_\infty  := \bigcap_{n\in \omega} p_n$ is in ${PT_H}$
  \ite 2 For all $n$: $p_n \weaker_{n+1}  p_\infty$. 
  \endent

\bigskip

\ppro\trimFact/Fact: 
\begindent
\ite 1  If  $\eta\in p \in {PT_H}$, then $ \trim p \eta \in {PT_H}$, and $p
   \weaker    \trim p \eta$.  (See \emNot/.) 
\ite 2 If $p \weaker q$ are conditions in ${PT_H}$, $\eta \in q$, then
   $\trim p \eta \weaker \trim q \eta $. 
\endent

\ppro\joinFact/Fact: If for all $ \eta\in \split_{n}(p)$, $q_\eta
   \stronger \trim p \eta$ is a condition in ${PT_H}$, then
\begindent
\ite 1  $q:=   {\bigcup\limits_{\eta\in \split_{n} (p)} q_\eta}$ is in
        ${PT_H}$,
\ite 2   $q \stronger_n p$ 
\ite 3   for all $\eta\in   \split_{n}(p)$,  $\trim q \eta {} = q_\eta$. 
   \endent

\ppro\finSplit/Fact: If $n\in \omega$, $p\in PT_H$, then $\split_n(p)$
is finite.  

Proof: This is the first time that we use the fact that each $H(n)$ is
a finite set:  Assume that the conclusion  is not true, so for some
$n$ and $p$, 
$\split_n(p)$ is        infinite. Then also 
   $$ T:= \{\eta\on k: \eta\in \split_n(p), k\le |\eta|\} \subseteq
    p$$ is infinite.  As $T$ is a finitely splitting tree, there has
    to be an infinite branch $b \subseteq T$. By \splitFac/(2), there
    is $\nu\in b \subseteq T$, $\height_p(\nu)> n$.  This is a
    contradiction to     $\splitFac/(1)$. 

\ppro\SacksFact/Fact:  $PT_H$ is strongly \oob/, i.e.: 

If $\name \alpha$ is a ${PT_H}$-name for an ordinal, $p\in {PT_H}$,
   $n\in \omega$, then there exists a finite set $A$ of ordinals and a
   condition $q\in {PT_H}$, $p \weaker_n q$, and $ q \forces \alpha
   \in A$.

   Proof: Let $C:= \split_{n}(p)$. $C$ is finite. 
   For each  node $\eta\in C$, let $q_\eta
   \stronger \trim p \eta$ be a condition such that for some ordinal
   $\alpha_\eta$ 
   $q_\eta \forces \name  \alpha = \alpha_\eta$.   Now let
   $$ q:= \bigcup_{ \eta\in C}     q_\eta
       \qquad \hbox {and}\qquad A:= \{\alpha_\eta: \eta\in C\} $$
    Since any extension of $q$ must be compatible with some
   $\trim q {\eta} $ (for some $\eta\in C$), $q \forces \name \alpha
   \in A$. 
\medskip

\ppro\PTcor/Corollary: ${PT_H}$ is proper (and indeed satisfies axiom
    A, so is  $\alpha$-proper for any
          $\alpha<\omega_1$) and \oob/.  Moreover, if $n\in
   \omega$, $p\in {PT_H}$, $\nt $ a name for a set of ordinals, then
   there exists a condition $q \stronger_n p $ such that 
\begindent
\ite 1 If $p\, \forces$``$\nt$ is finite'', then there is a finite set $A$
   such that $q \forces$``$\nt \subseteq A$''. 
\ite 2 If $p\, \forces$``$\nt $ is countable'', then there is a
countable set
   $A$ such that  $q \forces $``$\nt \subseteq A$''. 
\endent
      
   Proof: Use \SacksFact/ and \fusionFact/. 
   
Similarly to \SacksFact/ we can show: 

\ppro\RPprop/Fact: 
 Assume that $\name \alpha$ is a $RP$-name for an
   ordinal, $p\in RP$,  $n\in \omega$. 
   
   Then there exists a countable set $A$ of ordinals and a condition
   $q\in {PT_H}$, 
   $p \weaker_n q$, and $ q \forces \alpha \in A$. 

Proof: Same as the proof of \SacksFact/, except that now the set $C$
    and hence also the set $A$ may be countable.

\ppro\RPproper/Fact:  $RP$ is proper (and satifies axiom A). 
Proof: By \RPprop/ and \fusionFact/.

   \ppro\ptgenDef/Definition: Let $G \subseteq {PT_H}$ be a
          $V$-generic filter.  Then we let  $\name g$ be the
$PT_H$-name defined by 
   $$\name g :=  \bigcap_{p\in G} p $$ 
We may write $\name g_H$ or $\name g_{PT_H}$ for this name  $\name g$.
If $PT_H$ is 
the $\alpha$th iterand $Q_\alpha$ in an iteration, we write $\name
g_\alpha$ for $\name g_H$.

\ppro\ptgenFact/Fact: $\emptycondition_{PT_H}$ forces that
   \begindent
     \ite 0 $\name g$  is  a function with domain $\omega$,  
   \ite 1 $ \forall n \,  \name g(n) \in H(n)$. 
   \ite 2 For all $f\in V$, if\  $\forall n \,\, f(n) \in H(n)$ then 
		$\exists^\infty n \,\,f(n)=\name g(n)$. 
   \endent 
Furthermore,  for all $p\in PT_H$, 
\begindent
\ite 3 $p \forces {} $ ``$\{\name g\on n:n\in \omega\}$ is a branch
             through $p$.  
\ite 4 $p \forces \forall k \exists n \name g\on n \in \split_k(p)$
\endent

   Proof: (0) and (2) are  straightforward density  arguments.  (1)
   and (3) follow immedaitely from the definition of $\name g$.  (4)
   follows from (3) and \splitFac/(2), applied in $V^{PT_H}$.

\ppro\againRoths/Remark:   Since \SMZ/ is equivalent to

{\advance\leftskip by 2cm\advance\rightskip by 2 cm\relax 
for every $H: \omega \to \omega $, for every $F\in [\prod_n
   H(n)]^{\lless \c}$, there exists $f^*\in \fct$ such that
   for every $f\in F$ there are infinitely many $n$ satisfying
   $f(n)=f^*(n)$,
 \smallskip}
 \ptgenFact/(2) shows  that if we have $\c=\aleph_2$ and Martin's
   Axiom for the  
   forcing notions ${PT_H}$ 
   (for all $H$), then  we also have \SMZ/. (In fact the ``easy''
   implication ``$\repl$'' of this equivalence is sufficient.) 
  This can be achieved by
   a countable support iteration  of length $\aleph_2$ of forcing
   notions ${PT_H}$, with 
   the usual bookkeeping argument (as in  [\IteratedCohen]).

 We will show a stronger result in \addcor/: If $P:=P_{\omega_2}$
   is the limit of a countable support iteration $\<P_\alpha,
   Q_\alpha: \alpha<\omega_2>$, where ``many'' $Q_\alpha$ are of the
   form $PT_{H_\alpha}$ for some $H_\alpha$, then some bookkeeping
   argument can ensure that $V^{P}\thinks \Add(\smz)$.

\medskip

Since ${PT_H} $ is \oob/, it does not add Cohen reals.  The same is true
   for a countable support iteration of forcings of the form ${PT_H}$.
   However, in \bigproof/ we will have to consider a forcing iteration in
   which some forcing notions are of the form ${PT_H}$, but others do add
   an unbounded real.  To establish that even these iterations do not
   add Cohen reals, we will need the fact that the forcing notion 
   ${PT_H}$ preserves many ultrafilters.

\ppro\PTinterpret/Definition:  Let $Q$ be a forcing notion, $\nx$ a
    $Q$-name, $p\in Q$, $p \forces \nx \subseteq \omega$.   We say
    that $x^* \subseteq \omega$ is an {\it interpretation} of $\nx$
    (above $p$), if for all $n$ there is a condition $p_n \stronger p$
    such that $p_n \forces \nx \cut n = x^*\cut n$. 

\ppro\PTintFact/Fact: Assume $Q$, $p$, $\nx$ are as in \PTinterpret/.
    Then 
\begindent
\ite 1  There exists $x^* \subseteq \omega$ such that $x^*$ is an
    interpretation of $\nx $ above $p$. 
\ite 2  If $p \weaker p'$ and $x^* $ is an interpretation of $\nx$
    above $p'$, then $x^*$ is an interpretation of $\nx $ above $p$. 
\endent

\ppro\ultraLemma/Lemma:  ${PT_H}$ preserves P-points, i.e.: 
If\  $\U\in V$ is a P-point ultrafilter on
   $\omega$, then 
     $$ \forces_{{PT_H}} \hbox{``$\U$ generates an ultrafilter.''}$$

Proof: Assume that the conclusion is false. Then there is a
   ${PT_H}$-name $\nt$ for a subset of $\omega$ and a condition $p_0$
   such that  
   $$p_0\forces_{{PT_H}} \forall x \in \U: |x\cut \nt | = |(\omega - x)\cut
   \nt| = {\aleph_0} .$$

For each $p\in PT_H$ we choose a set $\tau(p)$ such that 
\begindent
\item{$\cdot$} $\tau(p)$ is an interpretation of $\ntau$ above $p$. 
\item{$\cdot$} If there is an interpretation of  $\ntau$ above $p$
    that is an element of $\U$, then $\tau(p)\in\U$. 
\endent

  Note that if $\tau(p)\in \U$ and $p \stronger p'$, then also
   $\tau({p'})\in \U$, since (by \PTintFact/(2)) we could have chosen
   $\tau({p'}):=\tau(p)$.  Hence either for all $p$ $\tau(p)\in \U$,
   or for some $p_1 \stronger p_0$, all $p \stronger p_1$,
   $\tau(p)\notin \U$.  In the second case we let $\name \sigma $ be a
   name for the complement of $\nt$, and let $\sigma(p)=\omega -
   \tau(p)$. Then $\sigma(p)\in \U$ for all $p \stronger p_1$. Also,
   $\sigma(p)$ is an interpretation of $\nsigma$ above $p$.

  So wlog for all $p \stronger p_1$, $\tau(p)\in \U$ for some $p_1\in
   {PT_H}$, $p_1 \stronger p_0 $. 

\medskip

  We will   show that there is a condition $q \stronger p_1$ and a
   set $a\in \U$ such that $q \forces a \subseteq \nt$. 

Recall that as $\U$ is a P-point, player NOTIN does not have a winning
strategy in the P-point game for $\U$ (see
\ppGame/).

We now define a strategy for player NOTIN.  On the side, player NOTIN
   will construct a fusion  sequence $\sqn p $ and  a sequence  $\sqn
   m$ of natural numbers. 

   $p_0$ is given. 

   Given $p_n$, we let 
   $$ A_n = \bigcap_{\eta\in \split_{n+1}(p_n)} \tau(\trim {p_n} \eta) $$ 
   This set is in $\U$.  Player IN  responds with a finite set
    $a_n \subseteq A_n$.  Let $m_n := 1 + \max(a_n)$.  For each $\eta
   \in \split_{n+1} (p_n)$ there is a condition $q_\eta \stronger
   \trim {p_n} \eta$ 
   forcing $\ntau \cut m_n = \tau(\trim {p_n} \eta)\cut m_n$, so in
   particular 
    $$ q_\eta \forces a_n \subseteq \nt  \cut m_n$$
   Let $p_{n+1} = \bigcup\limits_{\eta\in \split_{n+1}(p_n)} q_\eta$. 

   Then 
   $$p_{n+1} \stronger_{n+1}  p_n \hbox{\quad and\quad} p_{n+1}
   \forces  a_n \subseteq \nt  \leqno (*)$$ 

  This is a well-defined strategy for player NOTIN. As it is not a
   winning strategy, there is a play in which IN wins.  During this
   play, we have constructed a fusion sequence $\sqn p $. 
   Letting 
   $a:=\bigcup_n a_n$, $q:= \bigcap_n p_n$, we have that $a\in \U$,
   $p_0 \weaker q \in {PT_H}$ (by \fusionFact/), and $q \forces a
   \subseteq \nt$ (by $(*)$), a contradiction to our assumption.

\bigskip
 The following facts will be needed for the proof that if we iterate
   forcing notions ${PT_H}$ with carefully chosen functions $H$, then we
   will get a model where the ideal of strong measure zero sets is
   $\c$-additive. 
\bigskip

\def\keepit{u}
\def\throwout{v}

\ppro\ptwFact/Fact and Definition:   Assume $p\in {PT_H}$, $\keepit
   \subseteq 
   \omega$ is infinite, $\throwout = \omega - \keepit$.   Then we can
   define a stronger 
   condition $q$ by ``trimming'' $p$ at each node in
   $\split^\throwout(p)$. (See \splitDef/(4).) 
    Formally, let 
   ${\vec\imath}=\<i_\eta:\eta\in \split^\throwout(p)>$ be a 
   sequence 
   satisfying $i_\eta\in H(|\eta|) $ for all $\eta\in \split^\throwout
   (p)$. 

Then $$ p_{\vec \imath} := \{\eta\in p:\forall n\in\dom(\eta):
   \hbox{If $\eta\on n \in \split^\throwout (p)$, then
   $\eta(n)=i_{\eta\on   n}$}\}$$
 is a condition in ${PT_H}$%

Proof: Let $q:=p_{\vec\imath}$.  $q$ satisfies (A)--(B) of the
       definition \HPerfDef/ of ${PT_H}$. 
 The definition of $p_{\vec\imath}$
immediately implies: 
\begindent
\ite 1 If $\eta\in \split^{\throwout}(p)\cut q$, then
       $\succ_q(\eta)=\{i_\eta\}$. 
\ite 2 If $\eta\in \split^{\keepit}(p)\cut q$, then 
        $\succ_q(\eta)=\succ_p(\eta) = H(|\eta|)$. 
\ite 3 If $\eta\in q -\split(p)$, then $\eta\in p - \split(p)$, so
       $\succ_q(\eta) = \succ_p(\eta)$ is a singleton. 
\endent
Note that $\split(p) = \split^{\keepit}(p) \cup
  \split^{\throwout}(p)$, so (1)--(3) cover all possible cases for
  $\eta\in q$.

So $q$ also satisfies \HPerfDef/(C).  

\relax From (1)--(3) we can also conclude: 
\begindent
\ite 4 For all $\eta\in q$:   $\succ_q(\eta)  \not=\emptyset$. 
\endent

To show that $q\in PT_H$, we still have to check condition
  \HPerfDef/(D).  So let $\eta\in q$.    Since ${\keepit}$ is
  infinite, there is $k\in \keepit$, $k>|\eta|$.  By (4), there is an
  infinite branch $b \subseteq q $ containing $\eta$.  By
   \splitFac/(2) there is $\vu\in b$, $\height_p(\nu) = k$. Then $\eta
   \subseteq \nu$, and $\nu\in \split(q)$.
\medskip

\ppro\ptwtwo/Fact:  $p_{\vec\imath}\, \forces$``$\eta \subseteq
   {\ng} \logand \eta\in \split^\throwout(p) \ \limpl \
   {\ng}(|\eta|)=i_\eta$'' (where ${\ng}$ is a name for the generic branch
   defined in \ptgenFact/).

Proof: $p_{\ivec} \forces {\ng} \subseteq p_{\ivec}$ and
   $\succ_{p_{\vec\imath}}(\eta) = \{i_\eta\}$.

\medskip

To simplify notation,  we will now assume  that for all $n$,  $H(n)\in
\omega$.  (If  $H(n)$ are just  arbitrary finite sets as in
\HPerfDef/, then we could prove analogous statements, replacing $0$
and $1$ by any two elements $0_n \not= 1_n$ of $H(n)$.)

\ppro\splitsDef/Definition:  Let $\name f $ be a ${PT_H}$-name for a
   function from $\omega $ to $\omega$. 
We say that $\nf $ {\it splits} on $p,k$
   if for all $\eta\in \split_{ k}(p)$ there are $l$ and
   $j_1\not=j_0$ such that  
   $$ \eqalign{ \trim p {\eta\extend 0 } & \forces \nf(l)=j_0\cr 
		\trim p {\eta\extend 1 } & \forces \nf(l)=j_1\cr}$$

\ppro\splitRem/Remark: If $\nf $ splits on $p$, $k$, and $q \ge_{k+1}
   p$, then $\nf $ splits on $q$, $k$. 

(Proof: $\split_k(p)=\split_k(q)$, and for $\eta\in \split_k(p)$,
        $\trim p {\eta\extend i} \le \trim q {\eta\extend i} $.)

\ppro\splitOne/Lemma: If $p \forces \nf \notin V$, $k\in \omega$, then
   there is $q \stronger_{k+1} p$ such that  $\nf$ splits on $q$, $k$.

Proof:   For $\eta\in \split_k(p)$, $i\in \{0,1\}$ we let $\eta_i$ be the
   unique element of $\split_{k+1}(p)$ satisfying $\eta\extend i
   \subseteq \eta_i$. 

  By \rzeroone/, 
for each $\eta\in \split_k(p)$ we can
find conditions     $q_{\eta_0} \stronger \trim p {\eta_0} $,  
   $q_{\eta_1} \stronger \trim p {\eta_1} $ and integers $l_\eta$,
   $j_{\eta,0}\not= j_{\eta,1}$ such that $q_{\eta_0} \forces
   \nf(l_\eta)=j_0$,  $q_{\eta_1} \forces
   \nf(l_\eta)=j_1$.   If  $\nu\in \split_{k+1}(p) $ is not of the
   form $\eta_0$ or $\eta_1$ for any $\eta\in \split_k(p)$, then let
   $q_\vu:= \trim p \nu $.   

  By \joinFact/, $q:= \bigcup\limits_{\nu\in \split_{k+1}(p)} q_\vu$ is a
   condition, $q \stronger_{k+1} p$, and $q_\nu = \trim q \nu $ for all
   $\nu \in \split_{k+1}(p)$.

   We finish the proof of \splitOne/ by showing that $\nf$ splits on
   $q$, $k$:  Let $\eta\in \split_k(p)=\split_k(q)$.  Then $\trim q
   {\eta\extend 0 } = \trim q {\eta_0} = q_{\eta_0}$, so $\trim q
   {\eta\extend 0 } \forces \nf(l_\eta) = j_{\eta,0}$.  Similarly, 
   $\trim q    {\eta\extend 1 } \forces \nf(l_\eta) = j_{\eta,1}$.

\ppro\splitAll/Lemma: If $p \forces \nf \notin V$, then
   there is $q \stronger  p$, $\nf$ splits on $q$, $k$ for all $k$.  

Proof: By \splitOne/, \splitRem/ and \fusionFact/ (using a fusion
   argument). 
\bigskip

\ppro\funny/Lemma:  Assume $Q$ is a strongly \oob/ forcing
   notion. Let $\nf$ be a
$Q$-name for a function, $p$ a condition, $n 
   \in \omega$, $p \forces \nf\notin V$.  Then there exists a natural
   numer $k$ such that 
$$ \hbox{for    all $\eta\in \kzerooneseq k $ there is a
   condition $q \stronger_n p$, $q \forces \nf \notin
[\eta]$.}\leqno(*)$$ 

We will write $k_{p,n} $ or $k_{\name f, p,n}$ for the least such $k$.
Note that for any $k \ge k_{p,n}$, $(*)$ will also hold. 

Proof: Assume that this is false. So for some  $\nf$, $n_0$, $p_0$, 
 $$ \forall k\in \omega\, \exists \eta_k\in \kzerooneseq k : 
    \lnot( \exists q \stronger _{n_0} {p_0} \,\,\, q \forces \nf
   \notin[\eta_k])\leqno(\star)$$ 
  Let  $$ T:=
   \{\eta_k\on l: l\le k, k\in \omega\}$$
 $T$ is a finitely branching tree ($\subseteq \zerooneseq$) of
   infinite height, so it must have an infinite branch.  Let $f^*\in
   \zerooneseq $ be such that $\{f^*\on j:j\in \omega\} \subseteq T$.

   Since $f^*\in V$ but ${p_0}\forces_Q \nf\notin V$, there exists a name
   $\name m$ of a natural number such that ${p_0} \forces f^*\on \name m
   \not = \nf\on \name m$.  Let $q \stronger_{n_0} {p_0} $ be such that for
   some $m^*\in \omega$, $q \forces \name m < m^*$. 

   Claim:  For some $k$, $ q \forces \nf \notin [\eta_k]$.  This will
   contradict $(\star)$. 

   Proof of the claim:  We have $q \forces \nf\on m^* \not= f^*\on
   m^*$.  Since $f^*\on m^*\in T$, there is a $k \ge m^*$ such that
   $f^*\on m^* = \eta_k\on m^*$.   Hence $q \forces \nf\on m^* \not=
   f^*\on m^* = \eta_k\on m^*$, so $q \forces \nf \notin [\eta_k\on
   m^*]$.  But then also  $q \forces  \nf \notin [\eta_k]$.  

   This finishes the proof of the claim and hence of the lemma.

\ppro\noincoo/Lemma: Assume that $Q$ is a strongly \oob/ forcing
   notion, $\H$ is a dominating family in $V$, 
       and ${\bar\nu} = \<\nu^h:h\in \H>$ has index $\H$.  Then
       $$\forces_{Q} \bigcap_{h\in\H}\bigcup_{k\in \omega}[\vu^h(k)]
       \subseteq V$$

\begingroup\def\\{\eta_0\in \kzerooneseq \eta_0 ,
   \ldots, \eta_{n-1} \in \kzerooneseq l_{n-1} }
Proof:  Assume that for some condition $p$ and some $Q$-name $\nf$, 
$$p    \forces \nf\notin V\logand \nf \in \bigcap_{h\in \H}\bigcup_{n\in
   \omega} [\vu^h(n)].$$
   We will define a tree of conditions such
that along every branch we have a fusion sequence. 

Specifically, 
 we will define an infinite sequence
   $\<l_n:n\in \omega>$ of natural numbers, and for each $n$ a finite
   sequence 
    $$\<p(\eta_0, \ldots, \eta_{n-1}):\\> $$
 of conditions satisfying
\begindent
\ite 0 $p() = p$
\ite 1 For all $n$: $\forall \\: \,\, l_n\ge k_{p(\eta_0, \ldots,
   \eta_{{n-1}}),n} $. 
\ite 2 For all $n$: $\forall \\\, \forall \eta_n\in \kzerooneseq l_n $
\begindent
\ite a $p(\eta_0, \ldots, \eta_{n-1} ) \weaker_n p(\eta_0, \ldots,
   \eta_{n-1}, \eta_n)$. 
\ite b $ p(\eta_0, \ldots, \eta_{n-1}, \eta_n) \forces \nf
   \notin[\eta_n]$. 
\endent
 \endent
Given $p(\eta_0, \ldots, \eta_{n-1}) $ for all $\\$, we can find
$l_n$ satisfying condition (1).   The by the definition of
$k_{p(\eta_0, \ldots, \eta_{n-1}),n}$ we can (for all $\eta_n\in
\kzerooneseq l_n $) find 
$p(\eta_0, \ldots, \eta_{n-1}, \eta_n)$.

Now let $h\in \H$ be a function such that for all $n$, $h(n) > l_n$.
Define a sequence $\<\eta_n:n\in \omega > $ by $\eta_n:=\nu^h(n)\on
l_n$, and let $p_n:=p(\eta_0, \ldots, \eta_n)$. Then $p \weaker p_0
\weaker_0 p_1 \weaker_1 \cdots\,$, so there exists a condition $q$
extending all $p_n$.  So for all $n$, $q \forces \nf \notin [\eta_n]$.
But then also for all $n$, $q \forces \nf \notin [\nu^h(n)]$, a
contradiction.

\endgroup      
\medskip

Lemma  \noincoo/ will be needed later to show that if we iterate
focings of the form $PT_H$ together with random real forcing, after
$\omega_2$ many steps we obtain no strong measure zero sets of size
$\aleph_2$.    The proof (in \smallproof/) would be much easier if we
could omit ``strongly'' from the hypothesis of \noincoo/, i.e., if we 
could answer the following question positively: 

\ppro\ooproblem/Open Problem: 
Assume $\H \subseteq \fct$ is a dominating family (or even wlog
$\H=\fct$), and   $\bar\nu$ has index $\H$.  Let $Q$ be an \oob/
forcing notion.  Does this imply 
$$\forces_Q\bigcap_{h\in\H}\bigcup_n[\nu^h(n)]\subseteq V\ \hbox{?}$$

\ppro\hstarfact/Fact: 
\begingroup
\def\\#1 #2 {\bigcup_{#1\in \omega}[#2(#1)]}
 Assume $h^*:\omega\to \omega- \{0\}$,
$H^*(n)=\kzerooneseq h^*(n) $.  Let $\H \subseteq \fct$ be a
dominating family, and let $\bar\nu$ have index~$\H$.  Let $\name g^*$
be the name of the generic function added by $PT_{H^*}$. 

Then $$ \forces_{PT_{H^*}} \exists h\in \H  \,\,
 \\k \nu^h \subseteq \\n \ng^* $$

Proof:  Assume not, then there is a condition $p$ such that 
$$ p \forces \forall h\in \H\,\, 
 \\k \nu^h \not\subseteq \\n \ng^* \leqno\qquad(*)$$
Let $h\in\H$ be a function such that $\forall k\in \omega \forall
\eta\in \split_{2k+1}(p) \,\, h^*(|\eta|)\le h(k)$.   This function
$h$ will be a witness contradicting $(*)$. 

For $\eta\in \split_{2k+1}(p)$ let $i_\eta\in \succ_p(\eta) =
H^*(|\eta|) = \kzerooneseq h^*(|\eta|) $ be defined by 
$ i_\eta := \vu^h(k)\on h^*(|\eta|)$.   (Note that $\vu^h(k)\in
\kzerooneseq h(k) $ and $h(k) \ge h^*(|\eta|)$.)

Let $\ivec:=  \<i_\eta:\eta\in \split_{2k+1}(p), k\in \omega>$ and 
let $q:= p_{\ivec}$. 

Then $q \forces \forall n \forall k \,(\ng\on n \in \split_{2k+1}(p)
\limpl \ng(n) = i_{\ng\on n} \subseteq \vu^h(k))$ by \ptwtwo/.

 Since
also $q 
\forces \forall k \exists n\,\, \ng\on n \in \split_{2k+1}(p)$, we get $q
\forces \forall k \exists n \,[\nu^k(k)] \subseteq [\ng(n)]$. 
This contradicts~$(*)$. 

\endgroup

\bigskip\goodbreak\bigskip


{\bf \S 3\quad Two models of $\Add(\smz)$.}
\neusection

 Recall that $S^1_2 := \{{\delta}<{\omega_2}: cf({\delta})=
\omega_1\}$. 

\ppro\reflect/Lemma:  Let $\<P_\alpha, Q_\alpha:\alpha < \omega_2>$
    be an iteration of proper forcing noitions as in \ccLemma/,  
     $p\in 
    P_{\omega_2}$,  $\name A$ a 
    $P_{\omega_2}$-name.  If $p \forces{}$``$ \name A$ is a strong measure
    zero set,'' then there is a closed unbounded set $C \subseteq
    \omega_2$ and a sequence $\<\bar\nu_{\delta}: {\delta}\in C \cut
    S^1_2>$ such that each $\bar\nu_{\delta}$ is a $P_{\delta}$-name,
    and 
   $$ p\forces_{\omega_2} \hbox{$ \bar\nu_{\delta}$  has index $\fct\cut
    V_{\delta}$ and $\displaystyle\name A \subseteq \bigcap_{h\in
    \fct\cut V_{\delta}}\bigcup_{n\in \omega} [\vu^h(n)]$}$$

Proof:  Let $\name c $ be a $P_{\omega_2}$-name for a
function from $\omega_2$ to $\omega_2$ such that for all $\alpha<
\omega_2$, 
$$ \forces_{\omega_2} \forall h\in \fct\cut V_\alpha\, \exists \vu^h\in
V_{\name c(\alpha)}: \, \forall n\,\, \nu^h(n)\in \kzerooneseq h(n) \logand
\name A \subseteq \bigcup_n [\vu^h(n)]$$

(Why does $\name c$  exist?  Working in $V[G_{\omega_2}]$, note that
there are only $\aleph_1$ many functions in $\fct \cut V_\alpha$, and
for each such  $h$ there is a $\nu^h$ as required in
$\bigcup_{{\beta}<\omega_2}V_{\beta} $, by \nonew/.)

As $P_{\omega_2}$ satisfies the $\aleph_2$-cc, by \ccFact/(1) 
 we can 
    find a function $c\in V$ such that     $\forces_{\omega_2} \forall
    \alpha \name c(\alpha ) < c(\alpha)$.   Let  
$$ C:= \{{\delta}: \forall \alpha< {\delta}\,\, c(\alpha) <
{\delta}\}$$\def\nunameh{\name\nu^{\smash{\name h}}}%
The set $C$ is closed  unbounded.   In $V$, we can assign to each
$P_\alpha$-name $\name h$ (for $\alpha<{\delta}\in C$) a $P_{\delta}$-name
$\nunameh$ such that 
$$ \forces_{\omega_2} \forall n\, \nunameh(n)\in \kzerooneseq \name h(n) 
\logand 
\name A \subseteq \bigcup_n [\nunameh(n)]$$
Now in $V[G_{\delta}]$ we can choose  for each $h\in \fct$ an $\alpha<
{\delta}$ and a $P_\alpha$-name $\name h$ such that $h=\name
h[G_{\delta}]$. Then we let $\vu^h:= (\nunameh)[G_{\delta}]$.  Thus we
    found a sequence $\bar \nu = \<\nu^h:h\in V_{\delta}>
    \in V_{\delta}$  as required.

\ppro\smzisclub/Lemma: 

Assume $\<P_\alpha, Q_\alpha: \alpha < \omega_2>$ is a countable
   support iteration of proper forcing notions,  where for each
    ordinal ${\delta}\in S^1_2$ $\forces_{\delta}
    Q_{\delta}=PT_{H_{\delta}}$ for some $P_{\delta}$-name
    $H_{\delta}$.  We will write $g_{\delta}$ for the generic function
    added by $Q_{\delta}$. 

    Assume $\name\H$ is a name
   for a dominating family ($\subseteq \pre{\omega}{(\omega -
\{0\})}$) in $V_{\omega_2}$, and 
   $$\eqalign{ \forces_{\omega_2} \hbox{``For all $h\in \name \H$, }
S_{h}:=\{{\delta}<{\omega_2}: &
      cf({\delta})=\omega_1, Q_{\delta}={PT_{ H}}^{V_{{\delta}}}\}\cr
      &\hbox{ is    stationary\  (where $H(n)=\kzerooneseq h(n) $).''}} $$

Let $G_{\omega_2} \subseteq P_{\omega_2}$ be $V$-generic, then in
    $V[G_{\omega_2}]$, a set $A \subseteq \R$ is a strong measure zero
    set iff there is a closed unbounded set $C \subseteq {\omega_2}$
    such that for every ${\delta} \in C\cut S^1_2 $, $A \subseteq
    \bigcup_n [g_{\delta}(n)]$.

\medskip
Proof: First we prove the easy direction.  Assume that for some club
    $C$, for all ${\delta}\in C \cut S^1_2$, $A \subseteq
    \bigcup_n[g_{\delta}(n)]$.  Then for every $h\in V_{\omega_2}\cut
    \pre{\omega}{(\omega-\{0\})}$  there is a ${\delta} = {\delta}_h\in
    C\cut S_{h} \subseteq S^1_2$.  So  $Q_{\delta_{h}} =
    (PT_H)^{V_{\delta_h}} $, where $H(n) =  \kzerooneseq h(n) $.   Since
    $g_{{\delta}_h}(n)\in \kzerooneseq h(n) $, and $A \subseteq
    \bigcup [g_{\delta_h}(n)]$ for arbitrary $h$, $A$ is a strong
    measure zero set.

Now for the reverse implication:  In $V_{\omega_2}$, let  $A$ be a
    strong measure zero set.  By 
    the previous lemma, there is a club set $C \subseteq \omega_2$
    and a sequence $\<\bar\nu_{\delta}:{\delta}\in C\cut S^1_2>$  such
    that each  $\bar\nu\in V_{\delta}$ is a sequence with index $\fct
    \cut V_{\delta}$ and $V_{\omega_2} \thinks A \subseteq
    X_{\bar\nu_{\delta}}$.   By \hstarfact/ we have for all ${\delta}
    \in C\cut S^1_2$: 
     $$ V_{{\delta}+1} \thinks \exists h\in
    V_{\delta}\,\,\bigcup_n[\nu_{{\delta}}^h(n)] \subseteq  
\bigcup_n [g_{\delta}(n)]$$
So fix $h_0\in V_{\delta}$ witnessing this. 
This inclusion is absolute, so also  
 $$V_{{\omega_2}} \thinks \bigcup_n[\nu_{{\delta}}^{h_0}(n)] \subseteq 
\bigcup_n [g_{\delta}(n)]$$
Thus 
$$ V_{\omega_2} \thinks A \subseteq X_{\bar\nu_{{\delta}}}
\subseteq \bigcup_n [\nu^{h_0}_{\delta} (n)] \subseteq \bigcup_n
[g_{\delta}(n)]$$ and 
we are done.
 \bigskip

\ppro\addcor/Corollary: Assume $P_{\omega_2}$ is as above. 
Then $\forces_{P_{\omega_2}} \Add(\smz) $. 

Proof: Let $\<A_i:i\in \omega_1>$ be a family of strong measure zero
    sets in $V_{\omega_2}$.  To each $i$ we can associate a closed
    unbounded set $C_i$ as in \smzisclub/.  Let
    $\displaystyle C:=\bigcap_{i\in \omega_1}C_i$, then also $C$ is closed
    unbounded, and for all ${\delta}\in C \cut S^1_2$,
    $\displaystyle\bigcup_{i\in \omega_1} A_i \subseteq \bigcup_{n\in
    \omega} [g_{\delta}(n)]$.  Again by \smzisclub/,
     $\displaystyle\bigcup_{i\in \omega_1} A_i$ 
    is a strong measure zero set.

 \bigskip

Our first  goal is to show that $\SMZ/$ does not guarantee the
   existence 
   of a strong measure zero set of size $\c$. Clearly the model for
   this should satisfy $\d=\aleph_1$ (if $\c=\aleph_2$), so we will
   construct a countable support iteration of \oob/ forcing notions.

  \ppro\smallproof/Theorem: If \ZFC/ is consistent, then 

\centerline{\ZFC/ + $\c=\aleph_2$ +  $\smz = [\R]^{\lleq\aleph_1}$
     + no real is Cohen over $L$ }
\centerline{+ there is a generalized Sierpinski set}
  is consistent.

Proof: We will start with a ground model $V_0$ satisfying $V=L$.  Let
   $\H:= \pre{\omega}{(\omega-\{0\})}\cut L =
   \{h_\alpha:\alpha<\omega_1\}$, and let $H_\alpha(n)=\kzerooneseq
    h_\alpha(n) $. 

Let $\<S_\alpha:\alpha<\omega_1>$ be a family of disjoint stationary
   sets $\subseteq \{{\delta}<\omega_2: cf({\delta})=\omega_1\}$.

  Construct a countable support
   iteration $\<P_\alpha, Q_\alpha: \alpha < \omega_2>$  satisfying 
  \begindent
  \ite 1 For all even $\alpha<\omega_2$: 
  $$ \forces_{P_\alpha} \hbox{For some $h:\omega\to \omega-\{0\}$,
    letting $H(n)= \kzerooneseq h(n) $, $Q_\alpha={{PT_H}}$.}$$
  \ite 2 If ${\delta}\in S_\alpha$, then $\forces_{\delta}
   Q_{\delta}=PT_{H_\alpha}$.
\ite 3 For all  odd  $\alpha<\omega_2$: 
  $$ \forces_{P_\alpha} \hbox{  $Q_\alpha=$ random real forcing.}$$
   \endent

  By \preserveoo/ 
   (or as a consequence of  \strooFact/), 
   $P_{\omega_2}$ is \oob/, so $\forces_{\omega_2}$``$\H$ is a
   dominating family.''  By \ccLemma/(3) 
   and \ccFact/ 
    the assumptions 
   of \addcor/ are satisfied, so $\forces_{\omega_2} \Add(\smz)$.
   Also, $\forces_{\omega_2}$``$\c=\aleph_2$ and there are no Cohen
   reals over $L$.''   Letting $X$ be the set of  random reals added
at odd stages, $X$ is a generalized Sierpinski set:  Any null set
$H\in V_{\omega_2} $ is covered by some $G_{\delta}$ null set $H'$
that coded in some intermediate model.  As coboundedly many elements
of $X$ are random over this model, $|H\cut X| \le |H'\cut X| \le
\aleph_1$. 

To conclude the proof of \smallproof/, we have to show 
$$ V_{\omega_2} \thinks \hbox{``If $X \subseteq \R$ is of strong
   measure zero, then $|X|<\c$.''}$$


Since $\H$ is a dominating family, by \defdefFact/ it is enough  to
   show that in $V_{\omega_2}$ the     following holds: 
$$ \hbox{If ${\bar\nu}$  has index $\H$, then  $|X_{\bar\nu}| \le
       \aleph_1$.}$$

We will show: If ${\bar\nu} \in V_\alpha$  has index $\H$, 
   then $X_{\bar\nu}  \subseteq \fct \cut
   V_\alpha$.  (This is sufficient, by \nonew/.)

Assume to the contrary that $G_{\omega_2}$ is a generic filter,
   ${\bar\nu} \in V_\alpha$, and in
   $V[G_{\omega_2}]$ there is ${\delta}>\alpha$,  $f\in V_{\delta} -
   \bigcup_{{\gamma}<{\delta}}V_{\gamma}$, $f\in X_{\bar\nu}$.   So
   also $$V[G_{\delta}] \thinks f\in V_{\delta} -
   \bigcup_{{\gamma}<{\delta}}V_{\gamma} \hbox{ and } f\in X_{\bar\nu}
   $$

Let $\nf $ be a $P_{\delta}$-name, $\name{{\bar\nu}}$ a
   $P_\alpha$-name, and let $p\in P_{\delta}$
    be a condition forcing all this.  ${\delta}$ cannot
   be a successor ordinal, by \noincoo/.   So ${\delta}$ is a limit
   ordinal, 
   and $cf({\delta})$ must be $\omega$, otherwise we would have
   $\zerooneseq    \cut V_{\delta} 
   = \zerooneseq \cut \bigcup_{\gamma<\delta}V_{\gamma} $.

So we have reduced the problem to the following lemma: 

\ppro\noincLemma/Lemma:   Let \iter{} be a countable support iteration
   of forcings where each $Q_\alpha$ (for even $\alpha$) is of the
    form $PT_{H_\alpha}$ 
   for some ($P_\alpha$-name) $H_\alpha$, and $Q_\alpha$ is random
    real forcing for odd $\alpha$.     Let ${\delta} \le
   \varepsilon$ be 
   a limit ordinal of countable cofinality, and let $\nf $ be a
    $P_{\delta}$-name of a 
   function in $\zerooneseq$ such that $\forces_{\delta} \forall
   \alpha<{\delta} \,\, \nf \notin V_\alpha$.   

  Let $\H\in V_0$ be a dominating family of functions, and assume that
   $\bar \nu$ has index $\H$. 

  Then $\forces_{\delta} \nf \notin 
   \bigcap\limits_{h\in {\cal H}} 
		\bigcup\limits_{n\in \omega}
					   [\nu^h(n)]$. 
\medskip

For notational simplicity, we again assume that for all even
$\alpha$, $\forces_\alpha$``$H_\alpha: \omega \to \omega$ (rather than
$H_\alpha:\omega \to\finzerooneseq$).''

Before we prove this lemma, we need the following two definitions
(which make sense for any countable support iteration 
  $\<P_\alpha, Q_\alpha: \alpha <\omega_2>$).

\ppro\landDef/Definition and Fact: For $p\in P_{\varepsilon}$,
$\alpha<{\varepsilon}$, 
$p\on \alpha 
   \forces p(\alpha) \weaker \nr \in Q_\alpha$, we define $p\land
   \nr$ as follows:  $(p\land \nr)({\gamma})=p({\gamma})$ for
   ${\gamma}\not=\alpha$, and $(p\land \nr)(\alpha)= \nr$. 

Then $p\land \nr\in P_\varepsilon$, $p\land \nr \stronger p$, and
   $(p\land \nr)\on \alpha = p\on \alpha$, so in particular $p\on
   \alpha \forces p\land \nr \in P_{\varepsilon}/G_\alpha$.  Furthermore,
   $p\on \alpha \forces (p\land \nr)(\alpha) = \nr $. 

Also, if $p(\alpha)=\nr$, then $p\land \nr = p$.  

\ppro\landRdef/Definition and Fact:  If  $p\in P_\alpha$, $A$ a countable
   subset of ${\varepsilon}$, and $p \forces \nr \in
   P_{\varepsilon}/G_\alpha \logand \nr \stronger p\logand 
    \supp(\nr) \subseteq \alpha \cup A$, then  we define 
    $p\land \nr$ as follows:  

For ${\gamma}<\alpha$, $(p\land \nr)({\gamma})=p({\gamma})$.  For
   ${\gamma}\ge \alpha$ and ${\gamma}\in A$, $(p\land \name
   r)({\gamma})=r({\gamma})$.

Again, 
    $p\land \nr\in P_\varepsilon$, $p\land \nr \stronger p$, and
   $(p\land \nr)\on \alpha = p\on \alpha$, so in particular $p\on
   \alpha \forces p\land \nr  \in P_{\varepsilon}/G_\alpha$.

   Also, if  $p_1 \weaker
   p_2$, then $p_1\land \nr \weaker p_2 \land \nr$. 

\bigskip

\ppro\noincP/Proof of \noincLemma/: 
  $cf({\delta})=\omega $, so we
   can find an increasing sequence $\sqn {\delta} $  of even ordinals
    converging to 
   ${\delta}$.  Assume there is a condition $p$ forcing that $\nf \in
   \bigcap\limits_{h\in {\cal H}}
		\bigcup\limits_{n\in \omega}
					   [\nu^h(n)]$.

  We will define sequences 
  \vtop{\hbox{$\sqn p $,}\hbox{$\sqn F $,}\hbox{$\sqn {\ell} $,}%
 \hbox{$\<\name s_n:n\in \omega>$}%
 \hbox{$\<p^i_n:n\in\omega,i\in\{0,1\}>$,}}

such
   that the following hold:    For each $n$, $p_n$, $p^0_n$, $p^1_n$
   are  conditions in
   $P_{\delta}$, ${\delta}_n$ is an even ordinal $<{\delta}$, $F_n$ is
   a finite subset of ${\delta}_n$,  $\ell_n$ 
   is an integer, and $\name s_n$ is a $P_{\delta_n} $-name for an
   element of $\finSeq$.   (We let $p_0=p$, $F_0=\emptyset$,
   $\ell_0=0$, $p^1_0=p^0_0=p_0$, $s_0=\emptyset$). For all  $n>0$ we
   will have: 
\begindent
\ite 1 $p_{n-1}  \weaker_{F_n, n} p_{n} $. 
\ite 2 $F_n \subseteq {\delta}_n$, $F_{n-1}  \subseteq F_{n+1} $, 
    $\bigcup_k \supp(p_k )   \subseteq \bigcup_k    F_k$. 
\ite 3 ${\delta}_{n-1}  \in F_{n}$. 
\ite 4 $p_n\on {\delta_n} \forces \name s_n = \stem(p_n({\delta_n}))
		=\stem(p_{n-1}({\delta_n})) $
\ite 5 For $i\in \{0,1\}$, $p^i_n = p_n \land \trim p_n({\delta_n})
   {s_n\vxtend i}$.  
\ite 6 $p_n\on {\delta}_n \forces_{{\delta}_n}$``$ \exists l<\ell_n\,
   \exists j_0\not=j_1 \, \forall i\in \{0,1\}: 
        p^i_n \forces_{{\delta_n},{\delta}} \nf(l)=j_i$.''
\endent

Note that (5) implies: 
\begindent
\ite 5' $p_n\on {\delta}_n \forces p^i_n\in P_{\delta}/G_{\delta_n}$,
\endent
and (6) implies 
\begindent
\ite 6' 
For  all $\eta\in \kzerooneseq \ell_n $: $ p_n\on {\delta}_n \forces
   \exists i\in \{0,1\}: p^i_n 
   \forces_{{\delta_n}, {\delta}}  \nf\on \ell_n \not= \eta$
\endent
[Proof of (6) $\limpl $ (6'): In $V_{\delta_n}$, let $i\in \{0,1\}$ be
   such that $\eta(l)\not= j_i$, where $l$ is as in (6).]

Finally, let $q= \bigcup_n p_n$.  Then $q\on {\delta}_n \forces
   \stem(p_n({{\delta_n} })) = \stem(q({\delta_n}))= \name s_n$ by (1),
   (3) and (4) and \sameStem/. 
   Let $h^* \in
   {\cal H}$ be a function such that for all $n$, $\ell_n < h^*(n)$.
   So for all $n$, $\nu^{h^*}(n)\on \ell_n$ is a well-defined member
   of $\kzerooneseq \ell_n $.

  For each $n$, let $\name i_n$ be a $P_{{\delta}_n}$-name of an
   element of $\{0,1\}$ such that 
 \begindent
\ite  6'' $ p_n\on {\delta}_n \forces p^{\smash{\name
   i}_n}_{\vphantom{\tilde n}n} \forces
   \nf\on \ell_n \not= \nu^{h^*}(n) \on \ell_n$
\endent

Now define a condition $q'$ as follows:   For $\alpha\notin
   \{{\delta}_n:n\in \omega\}$, $q'(\alpha)=q(\alpha)$, and 
   $$ q'({\delta}_n) = \trim q({\delta})  {s_n\vxtend\!\name i_n }$$
   (This is a $P_{{\delta}_n}$-name.)  

 Claim:  $q' \stronger q \stronger p$ (this is clear)
   and $q'  \forces \nf  \notin \bigcup\limits_{n\in \omega}
   [\nu^{h^*}(n)]$.   

To prove this claim, let $G_{\delta}\subseteq P_{\delta}$ be a generic
   filter containing $q'$, and assume $f:=\nf[G_{\delta}]$ is in
   $[\nu^{h^*}(n)]$.  Let $i_n:= \name i_n [G_{\delta_n}]$.  Now
   $q\in G_{\delta}$ implies $p_n\in G_{\delta}$, so in particular
   $p_n\on {\delta_n} \in G_{\delta_n}$.  Note that
   $\stem(q({\delta}_n)) = \stem(p_n({\delta}_n))= s_n$, so $q'\in
   G_{\delta}\logand p_n\in G_{\delta}$ implies $p^{i_n}_n \in
   G_{\delta}$.  Also, by (6'') we    have $q' \forces \nf \notin
   \bigcup_n [s^{f^*}(n)]$, a contradiction.

\bigskip

This finishes the proof of \noincLemma/ modulo the construction of the
sequences $p_n$, $F_n$, etc. 

First we fix enumerations $\dom(r) =  \{\alpha^m_r: m\in
   \omega\}$ for all $r\in P_{\delta}$.   We will write $\alpha^m_n$
   for $\alpha^m_{p_n}$.

Assume $p_{n-1}$ is given.   Let 
   $F_n:= {\delta}_n \cut \left(F_{n-1} \cup \{\alpha_k^m: k<n, m<n\}
   \cup    \{{\delta}_{n-1}\}\right)$.   This will take care of (2)
   and (3). 

To define $p_n$, first work in $V[G_{\delta_n}]$, where $p_{n-1} \on
   {\delta_n} \in G_{\delta_n}$.

We let $s_n:= \stem(p_{n-1} ({\delta_n}))$.

We let $r_0:= p_{n-1} \land \trim p_{n-1}({\delta_n}) {s_n\vxtend 0}$, 
and    $r_1:= p_{n-1} \land \trim p_{n-1}({\delta_n}) {s_n\vxtend 1}$. 

By \rzeroone/,  
 we can find $l$, $j_0\not=j_1$ and
$r_0'$, $r_1'$ such    that       $r_i' \stronger r_i$, and $r_i'
\forces \nf(l)=j_i$.

We now define a condition $r\in P_{\delta} / G_{\delta_n}$ as
      follows: 
\begindent
\item{$\cdot$} $r \on {\delta_n} = p_{n-1} \on {\delta_n} $.  
\item{$\cdot$}  $r({\delta_n}) = r'_0({\delta_n}) \cup
         r'_1({\delta_n}) \cup 
         \bigcup \{ \trim p_{n-1}({\delta_n}) {s_n\vxtend i}: 
	 i\in \succ_{p_{n-1}({\delta_n})}(s_n)-\{0,1\}\} $.
\item{}  (So $\stem(r({\delta_n})) = s_n$.) 
\item{$\cdot$}  If
      ${\gamma} \in \supp(p_{n-1})\cup \supp(r'_0)\cup \supp(r'_1)$
   and ${\gamma}>{\delta_n}$,  we let
      $r ({\gamma})  $ be a $P_{\gamma}$-name such that 
    $$\eqalign{p_{n-1}\on {\delta_n} \forcess_{\delta_n} \,
     \forcess_{\delta_n,{\gamma}} &
             \hbox{``For $i$ in $\{0,1\}$: If $\name s_n\vxtend i
 		   \subseteq 
      g_{\delta_n}$, then $r({\gamma})=r'_i({\gamma})$,}\cr
      & \hbox{\hphantom{``}and if
      $g_{\delta_n}$ extends neither $\name s_n\vxtend 0$ nor $\name
		   s_n\vxtend 1$,}\cr
      & \hbox{\hphantom{``}then $r({\gamma}) = p_{n-1} ({\gamma})
         $}\cr} $$ 

\item{} (We write $g_{\delta_n}$ for $g_{Q_{\delta_n}}$, the branch
   added by    the forcing $Q_{\delta_n}$.)
\endent

 This is a condition in $P_{\delta}/G_{\delta_n}$.  
Note that we have the following: 
\begindent
\ite i $ \stem(p_n({\delta_n})) = s_n = 
   \stem(p_{n-1}({\delta_n})) $
\ite ii For $i\in \{0,1\}$, $r\land r'_i({\delta_n}) \ge r_i'$. 
\ite iii $r \stronger_{{\delta_n} {\delta}} p_{n-1}$. 
\endent

Coming back to $V$, we can find names $\nr$, \dots, such that the
      above is forced by $p_{n-1}\on {\delta_n}$.   Now let $\bar r$
      be a condition in $P_{\delta_n}$ satisfying the following: 

\begindent
\ite a $\bar r \stronger_{F_n,n} p_{n-1}\on {\delta_n}$. 
\ite b For some countable set $A \subseteq {\delta}$, $\bar r \forces
      \supp(\nr) \subseteq A$. 
\ite c For some $\ell_n\in \omega$, $\bar r \forces \name l < \ell_n$. 
\endent

We can find a condition $\bar r $ satisfying (a)--(c) by
\strooFact/.   

      Finally, let $p_n:= \bar r \land\nr$. 
      So $p_n\on {\delta_n}    = \bar r$. 

And let $p^i_n $ be defined by (5). 

\bigskip
Why does this work?   

First we check (1): $p_{n-1} \on {\delta_n} \weaker_{F_n,n} p_n\on
   {\delta_n}$ by 
   (a), and $p_{n-1} \weaker p_n$, because $p_n\on {\delta_n} \forces
   p_n = \bar r \land \nr \stronger \nr \stronger p_{n-1}$ (by
   (iii)).  So $p_{n-1} \weaker_{F_n,n} p_n$.

(2) and (3) are clear. 

Proof of (4): $p_n\on {\delta_n} \forces \stem(p_n({\delta_n})) = 
   \stem((\bar r \land \nr)({\delta_n})) =\stem( \nr({\delta_n}))=
   \name s_n$.  

(6):    Let $G_{\delta_n}$ be a generic filter containing $p_n\on
   {\delta_n}$.   Work in $V[G_{\delta_n}]$.  We write $r$ for $\name
   r[G_{\delta_n}]$, etc. 

We want to show $p^0_n \forces_{{\delta_n} {\delta}} \nf(l)=j_0$. 
($p^1_n \forces_{{\delta_n} {\delta}} \nf(l)=j_1$ is similar.)  As
   $r'_0 \forces_{{\delta_n} {\delta}} \nf(l)=j_0$, it is enough to
   see $p^0_n \stronger r'_0$.  

\medskip

First we note that $p^0_n \stronger p_n \stronger p_{n-1} $.  Also
   $p^0_n = p_n\land \trim p_n({\delta_n}) {s_n\vxtend 0} \stronger
   p_n = \bar r \land r   \stronger r$.  

Finally, $p^0_n({\delta_n}) = \trim r({\delta_n}) {s_n\vxtend 0} =
   r'_0({\delta_n})$.   

So $p^0_n = p^0_n \land r'_0({\delta_n}) \stronger r\land
   r'_0({\delta_n}) \stronger r'_0$, and we are done.

\bigskip\bigskip

Our next model will satisfy 
$$ (*)\qquad\qquad\SMZ/ \  + \ \d=\c=\aleph_2.$$
This in itself is very easy, as it is achieved by adding $\aleph_2 $
Cohen reals to $L$.  (Also Miller [\SomeProperties] showed that \SMZ/ +
   $\c=\aleph_2$ +
$\b=\aleph_1$ is consistent.)

Our result says that we can obtain a model for $(*)$ (and indeed,
satisfying $\Add(\smz)$) without adding Cohen reals.  In particular,
$(*)$ does not imply $\Cov(\meager)$.

\ppro\bigproof/Theorem:  Con(ZFC) implies 

\centerline{Con(ZFC 
 + $\c=\d=\aleph_2>\b$ + $\Add(\smz)$   +  no real is Cohen over $L$)}

Proof (sketch): We will build our model by a countable support
iteration of length 
$\omega_2$ where at each stage we either use a forcing of the form
${PT_H}$, or rational perfect set forcing.   A standard bookkeeping
argument ensures that the hypothesis of \addcor/ is satisfied, so
we get $\forces_{\omega_2} \Add(\smz)$.   Using rational perfect set
forcing on a cofinal set yields $\forces_{\omega_2} \d=\c=\aleph_2$.
Since all P-point ultrafilters from $V_0$ are preserved, no Cohen
reals are added. 

Proof (detailed version):   Let
$\{{\delta}<{\omega_2}:cf({\delta})=\omega_1\} \supseteq
\bigcup_{{\gamma}<\omega_2} S_{\gamma}$, where
$\<S_{\gamma}:{\gamma}<{\omega_2}>$ is a family of disjoint stationary
sets.  Let ${\Gamma} : \omega_2\times \omega_1 \to \omega_2$ be a
bijection.  We may assume that ${\delta}\in S_{\Gamma(\alpha,
{\beta})} \, \limpl \, {\delta}>\alpha$.

First we claim that there is a countable support iteration 
$\<P_\alpha, Q_\alpha: \alpha<\omega_2>$ and a sequence of names 
  $\< {\<{\name H_\alpha^{\beta}:
   \alpha<{\omega_2}}>:{\beta}<\omega_1}>$ such that

\begindent
\ite 1 For all $\alpha<{\omega_2}$, all ${\beta}<\omega_1$,
		$\name H_\alpha^{\beta}$ is a $P_\alpha$-name.
\ite 2 For all $\alpha<{\omega_2}$,
		$\forces_\alpha \{H^{\beta}_\alpha:{\beta}<\omega_1\}
		= \pre{\omega}{(\omega-\{0,1\})}$. 
\ite 3 For all ${\alpha}<{\omega_2}$:  If $\alpha\notin
		\bigcup_{{\gamma}<{\omega_2}}S_{\gamma}$, then
                $\forces_\alpha  Q_\alpha = RP$. 
\ite 4 For all $\alpha<{\omega_2}$, all ${\beta}<\omega_1$, all
		${\delta}\in S_{{\Gamma}(\alpha, {\beta})}$:
		$\forces_{\delta} Q_{\delta} =
		PT_{\name H_\alpha^{\beta}}$.  
\endent
Proof of the first claim:  By induction on $\alpha$ we can first
define $P_\alpha $, then $\<\name H^{\beta}_\alpha:{\beta}<\omega_1>$
(by \ccLemma/(1)), then $Q_\alpha$ (by (3) or (4),
   depending on whether $\alpha \in \bigcup_{{\gamma}<{\omega_2}}
   S_{\gamma}$ or not). 

\medskip
Our second claim is that letting $\name \H$ be a name for all
functions from $\omega $ to $\omega - \{0,1\}$ in $V[G_{\omega_2}]$,
the assumptions of \addcor/ are satisfied, namely: 
\begindent
\ite a $\Forces_{\omega_2}``$\forall H\in \name \H\, \exists {\gamma}<
{\omega_2}$ $S_{\gamma} \subseteq S_H$.''$
\ite b $\Forces_{\omega_2} ``$\forall {\gamma}<{\omega_2}$
$S_{\gamma}$ is stationary.''$ 
\endent
      (b) follows from \ccLemma/(3) and 
   \ccFact/, and 
   (a) follows    from  
    $$ \Forces_{\omega_2} ``For all $H\in \name\H$ there is
        $\alpha<{\omega_2}$ and ${\beta}<\omega_1$ such that $H=\name
        H_\alpha^{\beta}$.''$$
which in turn is a consequence of \nonew/. 

So by \addcor/, $V_{\omega_2} \thinks \Add(\smz)$.   

Let $G_{\omega_2} \subseteq P_{\omega_2}$ be a generic filter,
   $V_{\omega_2} = V[G_{\omega_2}]$.  

Again by
   \nonew/, every  $\H \subseteq \fct \cut V_{\omega_2}$ of size $\le
   \aleph_1$ is a subset
   of some $V_\alpha$, $\alpha<{\omega_2}$, so $\H$ cannot be a
   dominating family, as rational perfect set forcing $Q_{\alpha+1}$
   will introduce 
   a real not bounded by any function in $\H \subseteq V_\alpha
   \subseteq V_{\alpha+1}$.   Hence $\d=\c=\aleph_2$. 

Finally, any P-point ultrafilter from $V$ is generates an ultrafilter
   in $V_{\omega_2}$, so there are no Cohen reals over $V$.

This ends the proof of \bigproof/.

\bigskip\goodbreak\bigskip

\parindent=0.8cm

\def\key#1[#2]#3\par{\advance\refno by 1
               \edef#1{\number\refno}
\write\mgfile{\def\noexpand#1{#1}}
\medskip
\item{[#1]} #3}

\newcount\refno
\refno=0

REFERENCES.\nobreak
\advance\rightskip by 0cm plus 1 cm
\def\ams#1.{}

\def\paper#1,{{\it #1}, } \def\publ#1\publaddr#2\yr#3 {#1, #2 #3}
\def\inbook#1,{in:  #1, }  \def\nextreference#1\par{\bigskip \noindent#1}
\def\journall#1:#2(#3)#4--#5.{#1, {\bf #2} (#3), pp.#4--#5.}
\def\journalp#1:#2(#3)#4.{#1, vol~#2 (#3), p.~#4.}
\def\book#1,{{\it #1},} \def\vol#1 {Vol.~#1}
\def\pages#1--#2{pp.#1--#2}

\def\AnM{Annals of Mathematics}
\def\AnPAL{Annals of pure and applied logic}
\def\FM{Fundamenta Mathematicae}

\def\JSL{Journal of Symbolic logic}

\def\TAMS{Transactions of the American Mathematical Society}
\def\PAMS{Proceedings  of the American Mathematical Society}

\def\contemporary{\inbook Axiomatic Set Theory, Boulder, Co 1983, 143--159,
Contemporary Math 31, AMS, Providence RI 1984}

 \key\IteratedForcing[Ba:IF] J.~Baumgartner, \paper{Iterated forcing},
        \inbook{Surveys in set theory {\rm
        (A.~R.~D.~Mathias, editor)}}, London Mathematical
        Society Lecture Note Series, No.~8, Cambridge
        University Press, Cambridge, 1983.

 \key\BlassShelah[BS:PP] A.~Blass and S.~Shelah, \paper {There may be
        simple $P_{\aleph_1}$- and $P_{\aleph_2}$-points, and the
        Rudin-Keisler order may be downward directed}, 
        \journall\AnPAL:33(1987)213--243. \ams 03E55 (03E05).

 \key\Corazza[Co:GB] P.~Corazza, \paper{The generalized Borel
        Conjecture and strongly proper orders},
        \journall\TAMS:316(1989)115--140. . 

\key\ShelahsPreservation[GJ:SP] M.~Goldstern, H.~Judah,  On Shelah's
        preservation theorem, preprint. 

 \key\JSW[JSW:BC] H.~Judah, S~Shelah, H.~Woodin, \paper The Borel
         Conjecture, \AnPAL, to appear.

\key\smzRapid[Ju:SR] H.~Judah, \paper Strong measure zero sets and
        rapid filters,\journall\JSL:53(1988)393--402.

 \key\Kunen[Ku:ST] K.~Kunen, \book Set Theory: An Introduction to
         Independence Proofs,

 \key\Mapping[Mi:MS] A.~Miller, \paper Mapping a set of reals onto the
         reals, 
     \journall\JSL:48(1983)575--584. 

 \key\Rational[Mi:RP] A.~Miller, \paper Rational Perfect Set Forcing,
\contemporary.
\ams 03E35 (03E49, 06E10). 

 \key\SomeProperties[Mi:MC] A.~Miller, \paper Some properties of
        measure and category,  \journall\TAMS:266,1(1981)93--114. %
        \ams 03E35 (03E15, 28C15, 54A35, 54H05). 

 \key\Pawlikowski[Pa:TB] J.~Pawlikowski, \paper Power of transitive
         bases of measure and category,
         \journall\PAMS:93(1985)719--729.

 \key\ProprieteC[Ro:PC] Rothberger, \paper Sur des families
        indenombrables de suites de nombres naturels et les
        probl\`emes concernant la propriete\'e C,
        \journall{Proc.\ Cambr.\ Philos.\ Soc.}:37(1941)109--126. 

 \key\EigenschaftC[Ro:VC] Rothberger,  \paper Eine Versch\"arfung der
        Eigenschaft C, \journall\FM:30(1938)50--55.

 \key\ProperForcing[Sh:PF] S.~Shelah, \book Proper Forcing, Lecture Notes
        in Mathematics \vol942 , Springer Verlag.    

 \key\ProperImproper[Sh:PI] S.~Shelah, \book Proper and Improper
         Forcing, to appear in Lecture Notes      in Mathematics,
         Springer Verlag.    

 \key\IteratedCohen[ST:IC]  R.~Solovay and S.~Tennenbaum, \paper
         Iterated Cohen extensions and Souslin's  problem,
         \journall\AnM:94(1971)201--245.   

 \key\RealValued[So:RM] R.~Solovay, \paper Real valued measurable
         cardinals, \inbook{Axiomatic Set
        Theory}, Proc.\ Symp.\ Pure Math.\ {\bf 13 I} (D.~Scott, ed.),
        pp.397--428, AMS, Providence RI, 1971.

\bigskip\vfil

\halign{#\hfil&\tt\hfil#\cr
Martin GOLDSTERN &goldstrn@bimacs.cs.biu.ac.il\cr
Haim JUDAH & judah@bimacs.cs.biu.ac.il\cr
Department of Mathematics\span\omit\cr
Bar Ilan University\span\cr
52900 Ramat Gan, Israel\cr
Saharon SHELAH &shelah@shum.huji.ac.il\cr
Department of Mathematics\span\omit\cr
Givat Ram\cr
Hebrew University of Jerusalem\span\omit\cr
}
\bigskip
\bigskip
\

\eject
\bye